\documentclass{amsart}

\usepackage[headings]{fullpage}
\usepackage{amsmath} 
\usepackage{amssymb} 
\usepackage{amsopn}
\usepackage{amsthm}
\usepackage{amsfonts} 
\usepackage{mathrsfs}
\usepackage{mathtools}
\usepackage{stmaryrd}
\usepackage{manfnt}
\usepackage{mathdots}
\usepackage[OT2, T1]{fontenc}
\usepackage[dvips]{graphicx} \usepackage[usenames]{color}
\usepackage{xspace}
\usepackage[all]{xy}
\usepackage{dsfont}
\usepackage{ifthen}
\usepackage{marvosym}
\usepackage{wasysym}
\usepackage{marginnote}
\usepackage{float}
\usepackage{float}
\usepackage{booktabs}
\usepackage{soul}
\usepackage{colonequals}

\usepackage{subfigure}
\usepackage{tikz}
\usepackage{tikz-cd}
\usetikzlibrary{through,intersections,calc,angles,arrows}
\usetikzlibrary{shapes.geometric}
\usetikzlibrary{decorations.markings,calc}
\usepackage{pgfplots}
\pgfplotsset{compat=1.12}
\usepackage[pdfpagelabels, pdftex]{hyperref}
\hypersetup{
  pdftitle={},
  pdfauthor={},
  pdfsubject={},
  pdfkeywords={},
  colorlinks=true,    
  linkcolor=red,     
  citecolor=blue,     
  filecolor=blue,      
  urlcolor=blue,       
  breaklinks=true,
  bookmarksopen=true,
  bookmarksnumbered=true,
  pdfpagemode=UseOutlines,
  plainpages=false}
\DeclareRobustCommand{\gobblefive}[5]{}

\DeclareMathOperator{\rC}{C}

\DeclareMathOperator{\rF}{F}

\DeclareMathOperator{\rH}{H}

\DeclareMathOperator{\rR}{R}

\DeclareMathOperator{\rd}{d}

\newcommand{\bB}{{\mathbb B}}
\newcommand{\bC}{{\mathbb C}}

\newcommand{\bE}{{\mathbb E}}
\newcommand{\bF}{{\mathbb F}}

\newcommand{\bH}{{\mathbb H}}

\newcommand{\bQ}{{\mathbb Q}}

\newcommand{\bT}{{\mathbb T}}

\newcommand{\bZ}{{\mathbb Z}}

\newcommand{\cE}{{\mathscr E}}

\newcommand{\cO}{{\mathscr O}}

\newcommand{\dA}{{\mathcal A}}
\newcommand{\dB}{{\mathcal B}}

\newcommand{\dE}{{\mathcal E}}

\newcommand{\dG}{{\mathcal G}}

\newcommand{\dP}{{\mathcal P}}

\newcommand{\dT}{{\mathcal T}}

\newcommand{\dX}{{\mathcal X}}

\newcommand{\fU}{{\mathfrak U}}

\newcommand{\fX}{{\mathfrak X}}

\newcommand{\sB}{{\mathsf B}}

\newcommand{\sD}{{\mathsf D}}

\DeclareSymbolFont{cyrletters}{OT2}{wncyr}{m}{n}
\DeclareMathSymbol{\Sha}{\mathalpha}{cyrletters}{"58}

\DeclareMathOperator{\Aut}{Aut}

\DeclareMathOperator{\Cov}{Cov}
\DeclareMathOperator{\End}{End}

\DeclareMathOperator{\Hom}{Hom}

\newcommand{\Rep}{{\sf Rep}}

\DeclareMathOperator{\uHom}{\underline{Hom}}

\newcommand{\xyinj}{\ar@{^(->}}

\DeclareMathOperator{\GL}{GL}

\DeclareMathOperator{\Sp}{Sp}
 
\DeclareMathOperator{\rank}{rk}
\DeclareMathOperator{\Mod}{\mathbf{Mod}}

\DeclareMathOperator{\Spec}{Spec}

\DeclareMathOperator{\Ext}{Ext}

\DeclareMathOperator{\ind}{ind}

\def\10{{\overrightarrow{10}}}
\def\01{{\overrightarrow{01}}}

\newcommand{\alg}{{\rm alg}}
\newcommand{\an}{{\rm an}}

\newcommand{\cris}{{\rm cris}}

\newcommand{\pst}{{\rm pst}}
\newcommand{\deR}{{\rm dR}}

\newcommand{\et}{\mathrm{\acute{e}t}}
\newcommand{\proet}{\mathrm{pro\acute{e}t}}
\newcommand{\profet}{\mathrm{prof\acute{e}t}}

\newcommand{\nr}{{\rm nr}} 

\newcommand{\op}{{\rm op}}

\newcommand{\Zar}{{\rm Zar}}

\newtheorem{thm}{Theorem}[section]

\newtheorem{prop}[thm]{Proposition}
\newtheorem{lem}[thm]{Lemma}

\theoremstyle{definition}
\newtheorem{defi}[thm]{Definition}

\theoremstyle{remark}
\newtheorem{rmk}[thm]{Remark}

\newtheorem{ex}[thm]{Example}

\newenvironment{pro}[1][\proofname]{{\it{#1:}} }{\hfill $\square$}
\newenvironment{pro*}[1][\proofname]{{\it{#1:}} }{}
\newenvironment{pro**}[1][]{{\it{#1}} }{\hfill $\square$}

\numberwithin{equation}{section}
\usepackage{enum}
\newlist{enumer}{enumerate}{2}
\setlist[enumer]{label=(\roman*),align=left,labelindent=0pt,leftmargin=*,widest = (iii)}
\newlist{enumerar}{enumerate}{1}
\setlist[enumerar]{label=\arabic*.,align=left,labelindent=0pt,leftmargin=*,widest = 8.}
\newlist{enumera}{enumerate}{2}
\setlist[enumera]{label=(\arabic*),align=left,labelindent=0pt,leftmargin=*,widest = (8)}
\newlist{enumeral}{enumerate}{2}
\setlist[enumeral]{label=(\alph*),align=left,labelindent=0pt,leftmargin=*,widest = (m)}

\definecolor{shadecolor}{RGB}{186,238,186}
\definecolor{softlimegreen}{RGB}{186,238,186}
\definecolor{limegreen}{RGB}{208,243,208}
\definecolor{questioncolor}{RGB}{135, 173, 241} 
\definecolor{warningcolor}{RGB}{240,120,134}
\definecolor{verypaleyellow}{RGB}{255,255,194}

\newcommand{\Kbar}{\overline K}

\newcommand\llbrack{[\![}
\newcommand\rrbrack{]\!]}
\newcommand\un{\mathrm{un}}
\newcommand\gun{\mathrm{gun}}
\newcommand\Fil{\mathrm{Fil}}

\newcommand\gplike{\mathrm{gplike}}

\DeclareMathOperator{\DeR}{DR}
\newcommand{\FdeR}{\mathrm{FdR}}
\renewcommand{\st}{\mathrm{st}}

\renewcommand{\pro}{\operatorname{\mathrm{pro}}}
\newcommand{\hatotimes}{\operatorname{\hat\otimes}}

\newcommand{\ket}{\mathrm{k\acute{e}t}}
\newcommand{\proket}{\mathrm{prok\acute{e}t}}
\newcommand{\RH}{\mathrm{RH}}
\newcommand{\dRH}{\mathcal{RH}}

\renewcommand{\1}{\mathbf{1}}

\DeclareMathOperator{\fSet}{\mathbf{fSet}}
\renewcommand{\Cov}{\operatorname{\mathbf{Cov}}}
\renewcommand{\Vec}{\operatorname{\mathbf{Vec}}}
\renewcommand{\Rep}{\operatorname{\mathbf{Rep}}}
\DeclareMathOperator{\MIC}{\mathbf{MIC}}
\DeclareMathOperator{\FMIC}{\mathbf{FMIC}}
\DeclareMathOperator{\Loc}{\mathbf{Loc}}

\newcommand{\cts}{\mathrm{cts}}

\newcommand{\Rbar}{\overline R}

\newcommand{\Vbar}{\overline V}

\theoremstyle{theorem}
\newtheorem*{thm*}{Theorem}

\title{Local constancy of pro-unipotent Kummer maps}
\author{L.\ Alexander Betts}

\begin{document}

\setcounter{secnumdepth}{3}
\setcounter{tocdepth}{1}

\begin{abstract}
	It is a theorem of Kim--Tamagawa that the $\bQ_\ell$-pro-unipotent Kummer map associated to a smooth projective curve~$Y$ over a finite extension of~$\bQ_p$ is locally constant when~$\ell\neq p$. The present paper establishes two generalisations of this result. Firstly, we extend the Kim--Tamagawa Theorem to the case that~$Y$ is a smooth variety of any dimension. Secondly, we formulate and prove the analogue of the Kim--Tamagawa Theorem in the case~$\ell=p$, again in arbitrary dimension. In the course of proving the latter, we give a proof of an \'etale--de Rham comparison theorem for pro-unipotent fundamental groupoids using methods of Scholze and Diao--Lan--Liu--Zhu. This extends the comparison theorem proved by Vologodsky for certain truncations of the fundamental groupoids.
\end{abstract}
\maketitle

\tableofcontents

\section{Introduction}
\label{s:intro}

Let~$K$ be a field of characteristic~$0$ with algebraic closure~$\Kbar$, and let~$Y$ be a smooth geometrically connected variety over~$K$. For a prime number~$\ell$, let~$\Loc^\un(Y_{\Kbar},\bQ_\ell)$ denote the $\otimes$-category of unipotent $\bQ_\ell$-local systems on~$Y_{\Kbar,\et}$. This is a $\bQ_\ell$-linear Tannakian category, and any $\Kbar$-point $\bar x\in Y(\Kbar)$ determines a fibre functor $\omega^\et_{\bar x}\colon\Loc^\un(Y_{\Kbar},\bQ_\ell)\to\Vec(\bQ_\ell)$. The Tannaka groupoid of $\Loc^\un(Y_{\Kbar},\bQ_\ell)$ is called the \emph{$\bQ_\ell$-pro-unipotent \'etale fundamental groupoid} $\pi_1^\et(Y_{\Kbar};-,-)$ of~$Y_{\Kbar}$.

If we fix a choice of $K$-rational basepoint~$x_0$, then there is a natural continuous action of the absolute Galois group~$G_K$ on~$U\colonequals\pi_1^\et(Y_{\Kbar};\bar x_0,\bar x_0)$, as well as the Tannakian path-torsor $\pi_1^\et(Y_{\Kbar};\bar x_0,\bar y)$ for any~$y\in Y(K)$. If we pick a path $\gamma\in\pi_1^\et(Y_{\Kbar};\bar x_0,\bar y)(\bQ_\ell)$, then the map $\xi_\gamma\colon G_K\to U(\bQ_\ell)$ given by $\sigma\mapsto\gamma^{-1}\sigma(\gamma)$ is a non-abelian continuous cocycle. The class of~$\xi_\gamma$ is independent of the choice of~$\gamma$, and so there is a \emph{pro-unipotent Kummer map}
\[
j_\ell\colon Y(K) \to \rH^1(G_K,U(\bQ_\ell))
\]
sending~$y\in Y(K)$ to the class of~$\xi_\gamma$. This pro-unipotent Kummer map plays a central role in the Chabauty--Kim method \cite{minhyong:siegel,minhyong:selmer,minhyong-etal:bsd_conjecture}.

In some cases, the pro-unipotent Kummer map can be understood quite explicitly. From now on, let~$K$ be a finite extension of~$\bQ_p$ for some prime~$p$. When~$Y/K$ is a smooth projective variety with good reduction and $\ell\neq p$, then it follows from a standard specialisation argument that the pro-unipotent Kummer map~$j_\ell$ is constant. Without the good reduction and projectivity hypotheses, the pro-unipotent Kummer map~$j_\ell$ need not be constant, but is still well-behaved for curves by a theorem of Kim--Tamagawa.

\begin{thm*}[{\cite[Theorem~0.1]{minhyong-tamagawa}}]
	Suppose that~$K$ is a finite extension of~$\bQ_p$ for~$\ell\neq p$, and that~$Y/K$ is a smooth geometrically connected curve. Then the pro-unipotent Kummer map~$j_\ell$ is locally constant on~$Y(K)$ for the $p$-adic topology.
\end{thm*}

Our first (easy) result in this paper is a generalisation of the Kim--Tamagawa Theorem to arbitrary dimension.

\begin{thm}\label{thm:local_constancy_of_Kummer_maps_l-adic}
	Suppose that~$K$ is a finite extension of~$\bQ_p$ for~$\ell\neq p$, and that~$Y/K$ is a smooth geometrically connected variety. Then the pro-unipotent Kummer map~$j_\ell$ is locally constant on~$Y(K)$ for the $p$-adic topology.
\end{thm}


Our main interest in this paper, however, is to formulate and prove the analogue of these results in the case~$\ell=p$. Let us first observe that Theorem~\ref{thm:local_constancy_of_Kummer_maps_l-adic} does not hold verbatim in the case~$\ell=p$: for $Y/K$ any smooth projective curve with good reduction, the cohomology set~$\rH^1(G_K,U(\bQ_p))$ has a natural $\bQ_p$-scheme structure, and the Zariski-closure of the image of~$j_p$ is
\[
\rH^1_f(G_K,U(\bQ_p)) \colonequals \ker\left(\rH^1(G_K,U(\bQ_p))\to\rH^1(G_K,U(\sB_\cris))\right)
\]
where~$\sB_\cris$ is Fontaine's ring of crystalline periods \cite[Theorem~1]{minhyong:selmer}. So as soon as the genus of~$Y$ is $\geq1$, the map $j_p$ has infinite image and so cannot be locally constant.

The $p$-adic analogue of Theorem~\ref{thm:local_constancy_of_Kummer_maps_l-adic} which we will prove in this paper says that the pro-unipotent Kummer map~$j_p$ is locally constant ``modulo $\rH^1_e$''. To make this precise, consider the subset $\rH^1_g(G_K,U(\bQ_p))\subseteq \rH^1(G_K,U(\bQ_p))$ defined by
\[
\rH^1_g(G_K,U(\bQ_p)) \colonequals \ker\left(\rH^1(G_K,U(\bQ_p))\to\rH^1(G_K,U(\sB_\deR))\right) \,,
\]
and define an equivalence relation~$\sim_{\rH^1_e}$ on~$\rH^1(G_K,U(\bQ_p))$ by declaring $\xi\sim_{\rH^1_e}\eta$ just when~$\xi$ and~$\eta$ have the same image under the map
\[
\rH^1(G_K,U(\bQ_p)) \to \rH^1(G_K,U(\sB_\cris^{\varphi=1})) \,.
\]
We define
\[
\rH^1_{g/e}(G_K,U(\bQ_p)) \colonequals \rH^1_g(G_K,U(\bQ_p))/\sim_{\rH^1_e} \,.
\]
Our main theorem is as follows.

\begin{thm}\label{thm:local_constancy_of_Kummer_maps_p-adic}
	Suppose that~$K$ is a finite extension of~$\bQ_p$, and that~$Y/K$ is a smooth geometrically connected variety. Then the image of the pro-unipotent Kummer map~$j_p$ is contained in~$\rH^1_g(G_K,U(\bQ_p))$, and the composite
	\[
	Y(K) \xrightarrow{j_p} \rH^1_g(G_K,U(\bQ_p)) \to \rH^1_{g/e}(G_K,U(\bQ_p))
	\]
	is locally constant on~$Y(K)$ for the $p$-adic topology.
\end{thm}


\begin{rmk}
	We can even be quite precise about how locally constant the map in Theorem~\ref{thm:local_constancy_of_Kummer_maps_p-adic} is: if~$V\subseteq Y^\an$ is admissible open in the rigid analytification~$Y^\an$ of~$Y$, isomorphic to a closed polydisc, then the map is constant on~$V(K)$.
\end{rmk}

These two local constancy results play a central role in the theory developed in \cite{me:motivic_heights}, which relates these pro-unipotent Kummer maps to the classical theory of N\'eron--Tate heights on abelian varieties. Moreover, by generalising some of the foundational results supporting the Chabauty--Kim method, Theorems~\ref{thm:local_constancy_of_Kummer_maps_l-adic} and~\ref{thm:local_constancy_of_Kummer_maps_p-adic} constitute preliminary steps towards setting up the Chabauty--Kim method for higher-dimensional varieties and, more interestingly, for varieties with bad reduction.
\smallskip

Let us now say a little about the method of proof of Theorem~\ref{thm:local_constancy_of_Kummer_maps_p-adic}. For any $\Kbar$-point $\bar x\in Y(\Kbar)$, we have an associated \emph{universal pro-unipotent local system} ${}_{\bar x}\bar\bE^\et\in\pro{-}\Loc^\un(Y_{\Kbar},\bQ_p)$, namely the pro-object pro-representing the fibre functor~$\omega^\et_{\bar x}$. This universal pro-unipotent local system naturally carries the structure of a cocommutative coalgebra, induced by the tensor product on~$\Loc^\un(Y_{\Kbar},\bQ_p)$. When~$\bar x=x\in Y(K)$ is $K$-rational, ${}_{\bar x}\bar\bE^\et$ also carries a natural $G_K$-action, using which it descends to a pro-$\bQ_p$-local system~${}_x\bE^\et$ on~$Y_\et$. Theorem~\ref{thm:local_constancy_of_Kummer_maps_p-adic} boils down to proving the following.


\begin{enumerate}[label = \alph*) , ref = (\alph*)]
	\item\label{thmpart:pro_dR} For any two $K$-rational points $x,y\in Y(K)$, the fibre~${}_{\bar x}\bar\bE^\et_{\bar y}$ is pro-de Rham with respect to the natural action of~$G_K$.
	\item\label{thmpart:locally_constant} For fixed $x\in Y(K)$, the isomorphism class of~$\sD_\pst({}_{\bar x}\bar\bE^\et_{\bar y})$ (as a cocommutative coalgebra in the category of pro-$(\varphi,N,G_K)$-modules) is locally constant on~$Y(K)$ in the $p$-adic topology.
\end{enumerate}

The key ingredient to establish both of these facts is a comparison theorem between pro-unipotent \'etale and de Rham fundamental groupoids. Let~$\MIC^\un(Y,\cO_Y)$ denote the category of unipotent vector bundles with integrable connection on~$Y$. This is a $K$-linear Tannakian category, and any $K$-point $x\in Y(K)$ determines a fibre functor $\omega^\deR_x\colon\MIC^\un(Y,\cO_Y)\to\Vec(K)$. The \emph{de Rham fundamental groupoid} $\pi_1^\deR(Y;-,-)$ of~$Y$ is, by definition, the Tannaka groupoid of~$\MIC^\un(Y,\cO_Y)$ based at these fibre functors\footnote{It can of course happen that~$Y(K)=\emptyset$, in which case $\MIC^\un(Y,\cO_Y)$ is not obviously neutral and its fundamental groupoid as we've defined it is empty.}. As before, for any~$x\in Y(K)$ there is a \emph{universal pro-unipotent vector bundle with integrable connection} ${}_x\dE^\deR\in\pro{-}\MIC^\un(Y,\cO_Y)$ pro-representing~$\omega^\deR_x$, which carries the structure of a cocommutative coalgebra. We will prove the following.

\begin{thm}\label{thm:comparison}
	Suppose that~$K$ is a finite extension of~$\bQ_p$, and that~$Y/K$ is a smooth geometrically connected variety.
	\begin{enumerate}
		\item\label{thmpart:comparison_groupoid} For any $K$-rational points~$x,y\in Y(K)$, there is a canonical $G_K$-equivariant isomorphism
		\begin{equation}\label{eq:comparison_groupoid}
		\sB_\deR\otimes_{\bQ_p}\pi_1^\et(Y_{\Kbar};\bar x,\bar y) \cong \sB_\deR\otimes_K\pi_1^\deR(Y;x,y)
		\end{equation}
		of $\sB_\deR$-schemes.
		\item\label{thmpart:comparison_family} For any $K$-rational point~$x\in Y(K)$, the rigid analytification ${}_x\bE^{\et,\an}$ is a pro-de Rham local system on~$Y^\an_\et$ in the sense of \cite[Definition~8.3]{scholze:relative}, and there is a canonical isomorphism
		\begin{equation}\label{eq:comparison_family}
		\cO\bB_{\deR,Y^\an}\otimes_{\hat\bQ_p}{}_x\hat\bE^{\et,\an} \cong \cO\bB_{\deR,Y^\an}\otimes_{\cO_{Y^\an}}{}_x\dE^{\deR,\an}
		\end{equation}
		of pro-sheaves of $\cO\bB_{\deR,Y^\an}$-modules with integrable connection on the pro-\'etale site~$Y^\an_\proet$. (See \S\ref{ss:comparison} for precise definitions.)
	\end{enumerate}
	Moreover, the first of these isomorphisms is compatible with identities, composition and inversion in the fundamental groupoids, and the second is compatible with cocommutative coalgebra structures. Both isomorphisms are strictly compatible with Hodge filtrations (to be defined in \S\ref{ss:de_rham}).
\end{thm}

\begin{rmk}
	Theorem~\ref{thm:comparison}\eqref{thmpart:comparison_groupoid} generalises a number of \'etale--de Rham comparison theorems for fundamental groupoids available in the literature. Under various kinds of good reduction assumptions, similar comparison isomorphisms appear in \cite[Theorem~A]{vologodsky}, \cite[Theorem~1.8]{olsson}\cite[Corollary~8.13]{olsson_again}, and \cite[\S7]{hadian:motivic_pi_1}\cite[\S5]{faltings}. Beyond the good reduction case, a comparison isomorphism for semistable curves was constructed in \cite{andreatta-iovita-kim}.
	
	We remark that very general comparison theorems were proved in \cite{deglise-niziol} between different realisations of Cushman's motivic fundamental groupoid. It is presumably the case that the \'etale and de Rham realisations of the motivic fundamental groupoid can be identified with the \'etale and de Rham fundamental groupoids as defined via the Tannakian formalism, in which case Theorem~\ref{thm:comparison}\eqref{thmpart:comparison_groupoid} would be a special case of \cite[Corollary~4.19]{deglise-niziol}, but we will not prove this in this paper.
	
	
\end{rmk}


It is more or less clear that Theorem~\ref{thm:comparison}\eqref{thmpart:comparison_family} implies point~\ref{thmpart:pro_dR} above; what deserves a little more explanation is how it implies point~\ref{thmpart:locally_constant}. The extra ingredient required is theory developed by Shimizu \cite{shimizu:p-adic_monodromy}, which shows that for a horizontal de Rham local system~$\bE$ on a closed polydisc or spherical polyannulus~$V$, the isomorphism class of~$\sD_\pst(\bE_{\bar y}))$ as a $(\varphi,N,G_K)$-module is constant as $y$ ranges over~$K$-points of~$V$. Applying this theory to~${}_x\bE^\et$ then gives the local constancy required for Theorem~\ref{thm:local_constancy_of_Kummer_maps_p-adic}.

\begin{rmk}
	We do not actually need the full strength of Theorem~\ref{thm:comparison} for the results in this paper; in fact it just suffices to know that~${}_x\bE^{\et,\an}$ is a pro-de Rham local system. This can actually be proved rather quickly and indirectly, see e.g.\ \cite[Remark~3.2]{me-daniel:weight-monodromy}. However, Theorem~\ref{thm:comparison} implies something stronger about the pro-unipotent Kummer map~$j_p$: not only is it locally constant modulo $\rH^1_e$, but it is locally analytic and can be explicitly described in terms of iterated integrals on small discs, much as in the Chabauty--Kim method for curves. Although we don't expand on this here, we still prove the full statement of Theorem~\ref{thm:comparison} with an eye to future work.
\end{rmk}

\subsubsection*{Proof of the comparison theorem}

The perspective which we adopt in proving the Comparison Theorem~\ref{thm:comparison} is that comparison isomorphisms for fundamental groupoids should ultimately come from a comparison isomorphism on the level of Tannakian categories. We can make this precise using the logarithmic Riemann--Hilbert correspondence for rigid analytic varieties \cite{dllz:riemann-hilbert}.

Let us fix a smooth normal crossings compactification $(X,D)$ of~$Y$. For convenience, we replace the category $\Loc^\un(Y_{\Kbar},\bQ_p)$ by the equivalent category $\Loc^\un(X^\an_{\bC_K},\bQ_p)$ of unipotent Kummer \'etale $\bQ_p$-local systems on $X^\an_{\bC_K}$, where~$\bC_K$ is the completion of~$\Kbar$ and~$X^\an_{\bC_K}$ is given the log structure coming from the divisor~$D^\an_{\bC_K}$. Similarly, we replace~$\MIC^\un(Y,\cO_Y)$ by the equivalent category~$\MIC^\un(X^\an,\cO_{X^\an})$ of unipotent vector bundles with integrable logarithmic connection on~$X^\an$. Diao, Lan, Liu and Zhu define a certain ringed site~$(\dX,\cO_\dX)$ which plays the role of the ``base change'' of~$(X^\an_\an,\cO_{X^\an})$ from~$K$ to~$\sB_\deR$, as well as a category~$\MIC^\un(\dX,\cO_\dX)$ of unipotent $\cO_\dX$-vector bundles with integrable logarithmic connection. We show that~$\MIC^\un(\dX,\cO_\dX)$ is a Tannakian category over~$\sB_\deR$, and that every $K$-point~$x\in Y^\an(K)=Y(K)$ gives rise to a fibre functor $\omega^\RH_x\colon\MIC^\un(\dX,\cO_\dX)\to\Vec(\sB_\deR)$. Thus we have an associated fundamental groupoid $\pi_1^\RH(\dX;-,-)$, which we call the \emph{pro-unipotent Riemann--Hilbert fundamental groupoid}.

There are natural $\otimes$-functors
\[
\sB_\deR\hatotimes_K(-)\colon\MIC^\un(X^\an,\cO_{X^\an}) \to \MIC^\un(\dX,\cO_\dX) \hspace{0.4cm}\text{and}\hspace{0.4cm} \dRH\colon \Loc^\un(X^\an_{\bC_K},\bQ_p) \to \MIC^\un(\dX,\cO_\dX)
\]
($K$- and $\bQ_p$-linear, respectively), the latter of which is (essentially) the logarithmic Riemann--Hilbert functor \cite[Theorem~1.9]{dllz:riemann-hilbert}. The main input, coming from the theory of \cite{dllz:riemann-hilbert}, is that the induced maps
\[
\sB_\deR\otimes_K\Ext^1(\cO_{X^\an},\dE) \to \Ext^1(\cO_\dX,\sB_\deR\hatotimes_K\dE) \hspace{0.4cm}\text{and}\hspace{0.4cm} \sB_\deR\otimes_K\Ext^1(\bQ_p,\bar\bE) \to \Ext^1(\cO_\dX,\dRH(\bar\bE))
\]
on Yoneda Ext-groups are isomorphisms for all~$\dE\in\MIC^\un(X^\an,\cO_{X^\an})$ and all~$\bar\bE\in\Loc^\un(X^\an_{\bC_K},\bQ_p)$. This implies parts~\eqref{thmpart:comparison_groupoid} and~\eqref{thmpart:comparison_family} of Theorem~\ref{thm:comparison} in an essentially formal manner: the only part which requires any real care is checking compatibility with $G_K$-actions and Hodge filtrations.

\begin{rmk}
	One advantage of this approach -- constructing the comparison maps on the level of Tannakian categories -- is that the proof should easily generalise to more general kinds of fibre functors than just those arising from $K$-rational points. For example, it should be straightforward to give a version of Theorem~\ref{thm:comparison} where~$x$ and~$y$ are allowed to be $K$-rational tangential basepoints, but we do not pursue this here.
\end{rmk}


\subsection{Notation}

We fix the following notation for the rest of the paper. $p$ will be a prime number, and~$K$ will denote a finite extension of~$\bQ_p$. We fix an algebraic closure~$\Kbar$ of~$K$, and write~$\bC_K$ for its algebraic closure. We write~$G_K$ for the Galois group of~$\Kbar$ over~$K$, or equivalently the group of continuous automorphisms of~$\bC_K$ over~$K$.

We will denote by~$Y$ a smooth geometrically connected variety over~$K$, and choose a smooth compactification~$X$ of~$Y$ whose complementary divisor~$D$ has normal crossings. We view the rigid analytifications~$Y^\an$ and~$Y_{\bC_K}^\an$ as adic spaces in the usual way, and view~$X^\an$ and~$X_{\bC_K}^\an$ as log adic spaces in the sense of \cite[\S2.1]{dllz:foundations} by endowing them with the log structure coming from the analytification of~$D$ as in \cite[Example~2.3.16]{dllz:foundations}.

Rings and algebras are always commutative and unital; coalgebras are always counital. For a ring~$A$ we write~$\Mod(A)$ for the category of $A$-modules and~$\Vec(A)$ for the subcategory of finite projective $A$-modules.

A $\otimes$-category/functor/natural transformation means a symmetric monoidal category/functor/natural transformation. We will permit ourselves the usual imprecisions of not explicitly describing the compatibility isomorphisms on $\otimes$-categories or functors.

A Tannakian category~$\dT$ over a field~$F$ is an essentially small\footnote{This essential smallness requirement is omitted in some definitions but is necessary for a well-behaved theory of Tannakian categories, since a Tannakian category which is not essentially small cannot be equivalent to the category of representations of an affine group-scheme.} $F$-linear rigid abelian $\otimes$-category such that~$\End_\dT(\1)=F$ and $\dT$ admits a faithful exact $F$-linear $\otimes$-functor $\dT\to\Vec(A)$ for some $F$-algebra~$A$.

\subsection{Acknowledgements}

I am grateful to Jonathan Pridham for helpful discussions about comparison theorems for fundamental groups, and to Fabrizzio Andreatta and Ruochuan Liu for taking the time to answer my questions on various technical aspects of their papers. I am also grateful to Adam Morgan for introducing me to Breen's theory of cubical structures, and to Netan Dogra and Wies\l{}awa Nizio\l{} for drawing my attention to the comparison isomorphisms in~\cite{deglise-niziol}. Particular thanks are due to Minhyong Kim, without whose guidance and wisdom the present paper would not have been possible.

\smallskip

This article is adapted from part of the author's dissertation at the University of Oxford \cite{thesis}, and was substantially updated under the Simons Collaboration grant number 550031. The proof of Theorem~\ref{thm:local_constancy_of_Kummer_maps_p-adic} given in~\cite{thesis} is different to the one we give here: reducing to the case of semistable curves and using the comparison theory of \cite{andreatta-iovita-kim}.
\section{Local constancy of the $\bQ_\ell$-pro-unipotent Kummer map}
\label{s:main_theorem_l-adic}

We begin with the proof of the easier of our two main theorems: Theorem~\ref{thm:local_constancy_of_Kummer_maps_l-adic}. In fact, we will prove a stronger result, that if~$\pi_1^{(p')}(Y_{\Kbar},\bar x_0)$ denotes the maximal pro-prime-to-$p$ quotient of the profinite \'etale fundamental group of~$Y_{\Kbar}$, then the pro-$p'$ Kummer map
\[
j^{(p')}\colon Y(K) \to \rH^1(G_K,\pi_1^{(p')}(Y_{\Kbar},\bar x_0))
\]
(defined in exactly the same way as~$j_\ell$) is locally constant.

For this, let~$V$ be a rigid analytic space isomorphic to a closed polydisc and let~$V_{\bC_K}$ be its base-change to~$\bC_K$. We view~$V$ and~$V_{\bC_K}$ as adic spaces in the usual way \cite[(1.1.11)]{huber}. Let~$\Cov^{(p')}(V_{\bC_K})$ denote the category of prime-to-$p$ finite \'etale coverings of~$V_{\bC_K}$, i.e.\ finite \'etale maps $V'\to V_{\bC_K}$ such that~$V'$ is dominated by a finite \'etale covering of~$V_{\bC_K}$ which is Galois with group of order prime to~$p$. If~$\bar x\in V(\bC_K)$, then there is an associated fibre functor $\omega_{\bar x}^{(p')}\colon\Cov^{(p')}(V_{\bC_K})\to\fSet$ valued in the category of finite sets, and this makes $\Cov^{(p')}(V_{\bC_K})$. Its fundamental groupoid is the maximal pro-prime-to-$p$ quotient of the algebraic fundamental group of~$V_{\bC_K}$ in the sense of \cite{de_jong:rigid_pi_1}.

\begin{rmk}
	When we define~$\Cov^{(p')}(V_{\bC_K})$ as the category of prime-to-$p$ finite \'etale coverings of~$V_{\bC_K}$, it does not matter whether we view~$V_{\bC_K}$ as a rigid space over~$\bC_K$ or as an adic space by \cite[(1.1.11), Lemma~1.4.5(iv) and Proposition~1.7.11(i)]{huber}.
\end{rmk}

Now the absolute Galois group~$G_K$ acts on~$V_{\bC_K}$ in a natural way\footnote{This is where it is necessary that we view~$V_{\bC_K}$ as an adic space: the automorphisms of~$V_{\bC_K}$ given by elements of~$G_K$ are automorphisms of adic spaces, but are not morphisms of rigid spaces over~$\bC_K$.} and this induces an action on $\Cov^{(p')}(V_{\bC_K})$. If~$x\in V(K)$ is a $K$-rational point, then the fibre functor~$\omega^{(p')}_{\bar x}$ is $G_K$-invariant, and hence for any second $K$-rational point~$y$ there is an induced action of~$G_K$ on~$\pi_1^{(p')}(V_{\bC_K};\bar x,\bar y)$. Our main input is the following.

\begin{lem}\label{lem:trivial_paths_in_discs_l-adic}
	For any two $K$-rational points $x,y\in V(K)$ we have
	\[
	\pi_1^{(p')}(V_{\bC_K};\bar x,\bar y)^{G_K} \neq \emptyset \,.
	\]
	\begin{proof}
		If~$V$ is $0$-dimensional, then~$x=y$ and there is nothing to prove. If~$V$ is $1$-dimensional, then \cite[Theorem 6.3.2]{berkovich} and the discussion in \cite[\S4]{de_jong:rigid_pi_1} shows that every prime-to-$p$ finite \'etale Galois covering of~$V_{\bC_K}$ is trivial, i.e.\ is a disjoint union of copies of $V_{\bC_K}$. This implies that~$\pi_1^{(p')}(V_{\bC_K};-,-)$ is the trivial groupoid: $\pi_1^{(p')}(V_{\bC_K};\bar x,\bar y)$ consists of a single point, which is perforce $G_K$-fixed.
		
		In general, we would like to run the same argument, but unfortunately the author does not know a proof of the triviality of the pro-prime-to-$p$ fundamental groupoid of a polydisc. Instead, we observe that there is a morphism $V_1\to V$ from the $1$-dimensional disc~$V_1$ into~$V$ which takes some $K$-points $x_1,y_1\in V_1(K)$ to $x,y\in V(K)$. The induced map
		\[
		\pi_1^{(p')}(V_{1,\bC_K};\bar x_1,\bar y_1) \to \pi_1^{(p')}(V_{\bC_K};\bar x,\bar y)
		\]
		is then~$G_K$-equivariant, so we are done by reducing to the $1$-dimensional case.
	\end{proof}
\end{lem}

As an immediate consequence, we obtain a proof of Theorem~\ref{thm:local_constancy_of_Kummer_maps_l-adic}. Since~$Y$ is smooth, for any~$x\in Y(K)$ we may choose an admissible open $V\subseteq Y^\an$ containing~$x$ and isomorphic to a closed polydisc. If $\Cov^{(p')}(Y_{\Kbar})$ denotes the category of prime-to-$p$ finite \'etale coverings of~$Y_{\Kbar}$, then we obtain by base-change, analytification and pullback a restriction functor
\[
\Cov^{(p')}(Y_{\Kbar}) \to \Cov^{(p')}(V_{\bC_K})
\]
inducing for every~$\bar y\in V(\bC_K)\subseteq Y(\bC_K)$ a map
\[
\pi_1^{(p')}(V_{\bC_K};\bar x,\bar y) \to \pi_1^{(p')}(Y_{\Kbar};\bar x,\bar y) \,.
\]
When~$\bar y=y\in V(K)\subseteq Y(K)$ is $K$-rational, this map is $G_K$-equivariant, so we deduce from Lemma~\ref{lem:trivial_paths_in_discs_l-adic} that $\pi_1^{(p')}(Y_{\Kbar};\bar x,\bar y)^{G_K}\neq\emptyset$.

In particular, if we choose a $G_K$-invariant path $\gamma\in\pi_1^{(p')}(Y_{\Kbar};\bar x,\bar y)^{G_K}$, then composition with~$\gamma$ gives a $G_K$-equivariant isomorphism
\[
\pi_1^{(p')}(Y_{\Kbar};\bar x_0,\bar x) \xrightarrow\sim \pi_1^{(p')}(Y_{\Kbar};\bar x_0,\bar y)
\]
of $\pi_1^{(p')}(Y_{\Kbar},\bar x_0)$-torsors. This means that~$x$ and~$y$ have the same image under the pro-prime-to-$p$ Kummer map~$j^{(p')}$, i.e.\ $j^{(p')}$ is constant on the $p$-adic neighbourhood $V(K)$ of~$x\in Y(K)$. The corresponding statement for the $\bQ_\ell$-pro-unipotent Kummer map follows immediately.\qed
\section{The unipotent Tannakian formalism}\label{s:tannakian}

In this section, we collect a few facts about Tannakian categories that will be used in the sequel. Throughout this section, we fix a field~$F$ and an $F$-linear Tannakian category~$\dT$. We denote fibre functors on~$\dT$ by~$\omega_x$, $\omega_y$ etc., and adopt the shorthand~$E_x\colonequals\omega_x(E)$.

We write~$\ind{-}\dT$ and~$\pro{-}\dT$ for the categories of ind- and pro-objects in~$\dT$, respectively. We denote the induced tensor products on $\pro{-}\dT$ and~$\ind{-}\dT$ by~$\hatotimes$ and~$\otimes$, respectively, and denote by~$(-)^\vee\colon\pro{-}\dT\xrightarrow\sim(\ind{-}\dT)^\op$ the functor induced by duality on~$\dT$. For a fibre functor~$\omega_x\colon\dT\to\Vec(F)$, we also denote by~$\omega_x$ the evident extensions
\begin{align*}
	\ind{-}\dT &\to \ind{-}\Vec(F) = \Mod(F) \\
	\pro{-}\dT &\to \pro{-}\Vec(F) = \Mod(F)^\wedge \,,
\end{align*}
where~$\Mod(F)^\wedge$ denotes the category of topological $F$-modules which are isomorphic to a product of copies of~$F$ with the discrete topology.

\subsection{Universal objects}\label{ss:tannakian_universal_objects}

Any fibre functor~$\omega_x$ on~$\dT$ is pro-representable by~\cite[I.8.10.14]{sga4}, meaning that there is a pro-object~${}_xE^\dT\in\pro{-}\dT$ and an element $e_x^\dT\in{}_xE^\dT_x$ such that the map
\[
\Hom_{\pro{-}\dT}({}_xE^\dT,E) \to E_x
\]
given by evaluation at~$e_x^\dT$ is bijective for all~$E\in\dT$. We call ${}_xE^\dT$ the \emph{universal object} of~$\dT$. There are comultiplication and counit maps
\[
\Delta\colon {}_xE^\dT \to {}_xE^\dT\hatotimes{}_xE^\dT
\hspace{0.4cm}\text{and}\hspace{0.4cm}
\eta\colon {}_xE^\dT \to \1
\]
given by $\Delta(e_x^\dT)=e_x^\dT\hatotimes e_x^\dT$ and $\eta(e_x^\dT)=1$, respectively. These make~${}_xE^\dT$ into a cocommutative coalgebra object in $\pro{-}\dT$.

If now~$\omega_y$ is a second fibre functor, possibly equal to~$\omega_x$, then there is a canonical bijection
\begin{equation}\label{eq:tannakian_yoneda}
\Hom(\omega_x,\omega_y) \xrightarrow\sim {}_xE^\dT_y
\end{equation}
courtesy of the Yoneda Lemma, taking a natural transformation~$\alpha\colon\omega_x\to\omega_y$ (not necessarily $\otimes$-natural) to~$\alpha_{{}_xE^\dT}(e_x^\dT)\in{}_xE^\dT_y$. In fact,~\eqref{eq:tannakian_yoneda} is an isomorphism of pro-finite-dimensional vector spaces where the pro-structure on~$\Hom(\omega_x,\omega_y)$ is the one coming from its description as an end.

It follows that there are natural morphisms
\[
\mu_{x,y,z}\colon {}_yE^\dT_z\hatotimes{}_xE^\dT_y \to {}_xE^\dT_z
\hspace{0.4cm}\text{and}\hspace{0.4cm}
S_{x,y} \colon {}_xE^\dT_y \xrightarrow\sim {}_yE^\dT_x \,,
\]
the former corresponding to composition of natural transformations and the latter sending a natural transformation~$\alpha\colon\omega_x\to\omega_y$ to the natural transformation~$S_{x,y}(\alpha)\colon\omega_y\to\omega_x$ with components~$S_{x,y}(\alpha)_E = \alpha_{E^\vee}^\vee$.

We will see shortly that~$\mu_{x,y,z}$ and~$S_{x,y}$ are morphisms of cocommutative coalgebras, and that they make the~${}_xE^\dT_y$ into a groupoid in the category of cocommutative coalgebras.

\subsubsection{Description of the Tannakian fundamental groupoid}

The Tannakian fundamental groupoid of~$\dT$ can be described in terms of the universal objects~${}_xE^\dT$.

\begin{prop}\label{prop:tannakian_groupoid_via_universal_objects}
	For any two fibre functors~$\omega_x,\omega_y$ there is a canonical isomorphism
	\begin{equation}\label{eq:tannakian_groupoid_via_universal_objects}
	\pi_1(\dT;\omega_x,\omega_y) \cong \Spec({}_xE^{\dT,\vee}_y)
	\end{equation}
	of affine $F$-schemes. Moreover, under this identification, the element~$e_x^\dT\in{}_xE^\dT_x$ corresponds to the identity in~$\pi_1(\dT,\omega_x)(F)$, and the maps~$\mu_{x,y,z}$ and~$S_{x,y}$ correspond to the maps on duals of affine rings induced by composition and inversion in the fundamental groupoid of~$\dT$.
	\begin{proof}
		For any~$E_1,E_2\in\dT$ there is a map
		\[
		\Delta_{E_1,E_2}\colon\Hom(\omega_x,\omega_y)\to\Hom(E_{1,x},E_{1,y})\otimes\Hom(E_{2,x},E_{2,y})=\Hom(E_{1,x}\otimes E_{2,x},E_{1,y}\otimes E_{2,y}) \,,
		\]
		given by $\alpha\mapsto\alpha_{E_1\otimes E_2}$. Taking the inverse limit (or more accurately the end) over all~$E_1,E_2$ gives a map
		\[
		\Delta\colon\Hom(\omega_x,\omega_y) \to \Hom(\omega_x,\omega_y)\hatotimes\Hom(\omega_x,\omega_y)
		\]
		in the category of pro-finite dimensional $F$-vector spaces. It is easy to check that~$\Delta$ corresponds under the Yoneda isomorphism~\eqref{eq:tannakian_yoneda} to the comultiplication on~${}_xE^\dT_y$, and that the counit on~${}_xE^\dT_y$ corresponds to the map $\eta\colon\Hom(\omega_x,\omega_y)\to F$ given by the action of natural transformations on the unit object~$\1$.
		
		Hence the Yoneda isomorphism~\eqref{eq:tannakian_yoneda} is an isomorphism of cocommutative coalgebras in the category of pro-finite dimensional $F$-vector spaces. In particular, for any $F$-algebra~$\Lambda$, there is an induced bijection
		\[
		\Hom(\omega_x,\omega_y)_\Lambda^\gplike \xrightarrow\sim {}_xE^{\dT,\gplike}_{y,\Lambda}
		\]
		between the grouplike\footnote{An element of a coalgebra is called \emph{grouplike} just when it satisfies $\Delta(\gamma)=\gamma\otimes\gamma$ and $\Delta(\gamma)=1$.} elements of the base-changes to~$\Lambda$. But it follows from the above description that the grouplike elements of~$\Hom(\omega_x,\omega_y)_\Lambda=\Hom(\Lambda\otimes_F\omega_x,\Lambda\otimes_F\omega_y)$ can be canonically identified with the \emph{$\otimes$-natural} transformations $\Lambda\otimes_F\omega_x\to\Lambda\otimes_F\omega_y$, while the grouplike elements of ${}_xE^\dT_{y,\Lambda}=\Hom_F({}_xE^{\dT,\vee}_y,\Lambda)$ can be canonically identified with the \emph{$F$-algebra} maps ${}_xE^{\dT,\vee}_y\to\Lambda$. Since every $\otimes$-natural transformation $\Lambda\otimes_F\omega_x\to\Lambda\otimes_F\omega_y$ is a $\otimes$-natural isomorphism \cite[Proposition~1.13]{milne}, this shows that~$\Spec({}_xE^{\dT,\vee}_y)$ represents the functor of $\otimes$-natural transformations $\omega_x\to\omega_y$, and the result follows.
	\end{proof}
\end{prop}

\subsection{Unipotent Tannakian categories}\label{ss:tannakian_unipotent}


\begin{defi}
	We say that an object~$E$ in an abelian $\otimes$-category~$\dA=(\dA,\otimes,\1)$ is \emph{unipotent} just when it has a finite filtration
	\[
	0=\Fil_{-1}E\leq\Fil_0E\leq\Fil_1E\dots\leq\Fil_nE
	\]
	whose graded pieces are isomorphic to finite direct sums of the unit object~$\1$. The smallest possible value of~$n$ is called the \emph{unipotency class} of~$E$. We say that an $F$-linear Tannakian category~$\dT$ is \emph{unipotent} just when every object is unipotent; we say that~$\dT$ is \emph{finitely generated unipotent} if additionally $\Ext^1_\dT(\1,\1)$ is finite-dimensional over~$F$. This implies that in fact~$\Ext^1_\dT(E_1,E_2)$ is finite-dimensional for all~$E_1,E_2\in\dT$.
\end{defi}

\begin{rmk}
	If~$\dT$ is a Tannakian category and~$\omega_x$ a neutral fibre functor, then $\dT$ is unipotent (resp.\ finitely generated unipotent) if and only if its Tannaka group~$\pi_1(\dT,\omega_x)$ is pro-unipotent (resp.\ finitely generated pro-unipotent).
\end{rmk}

We note the following criteria to be a unipotent Tannakian category.

\begin{lem}\label{lem:tannakian_criteria}
	Let~$(\dA,\otimes,\1)$ be an $F$-linear abelian closed $\otimes$-category, and let~$\dT$ be the full subcategory consisting of unipotent objects.
	\begin{enumerate}
		\item\label{condn:tannakian_endos_of_unit} If~$\End(\1)=F$, then~$\dT$ is an abelian subcategory of~$\dA$.
		\item\label{condn:tannakian_closed_tensor} Suppose that for every exact sequence
		\[
		0 \to M_1 \to M_2 \to \1 \to 0
		\]
		in~$\dA$ and object~$M\in\dA$, the sequences
		\[
		0 \to M_1\otimes M \to M_2\otimes M \to M \to 0
		\hspace{0.4cm}\text{and}\hspace{0.4cm}
		0 \to M \to \uHom(M_2,M) \to \uHom(M_1,M) \to 0
		\]
		are again exact. Then~$\dT$ is a closed $\otimes$-subcategory of~$\dA$. Moreover, every object of~$\dT$ is dualisable and hence flat in~$\dA$.
		\item\label{condn:tannakian_lax_fibre_functor} Suppose that~\eqref{condn:tannakian_endos_of_unit} and~\eqref{condn:tannakian_closed_tensor} hold. Suppose that~$\omega\colon\dA\to\dB$ is an exact $F$-linear lax $\otimes$-functor valued in an $F$-linear abelian closed $\otimes$-category~$\dB$. Assume that~$\omega(\1)=\1_{\dB}$ is the unit object of~$\dB$. Then~$\omega$ restricts to an exact $\otimes$-functor $\dT\to\dB$. If~$\dB$ is not the zero category, this induced functor is moreover faithful.
	\end{enumerate}
	In particular, if~\eqref{condn:tannakian_endos_of_unit} and~\eqref{condn:tannakian_closed_tensor} hold, if there exists a functor~$\omega_x\colon\dA\to\Mod(A)$ as in~\eqref{condn:tannakian_lax_fibre_functor} for some non-zero $F$-algebra~$A$, and if $\Ext^1_{\dA}(\1,\1)$ is a set, then~$\dT$ is a unipotent Tannakian category over~$F$ and the restriction of~$\omega_x$ is a fibre functor.
	\begin{proof}
		\eqref{condn:tannakian_endos_of_unit} is \cite[Proposition~1.2.1]{shiho_1}.
		
		For~\eqref{condn:tannakian_closed_tensor}, let~$\dP$ denote the subcategory of~$\dA$ consisting of those objects~$P$ such for that the functors~$P\otimes(-)$ and $\uHom(P,-)$ are exact, and for every exact sequence
		\[
		0 \to M_1 \to M_2 \to P \to 0
		\]
		and every object~$M\in\dA$, the sequences
		\[
		0 \to M_1\otimes M \to M_2\otimes M \to P\otimes M \to 0
		\hspace{0.4cm}\text{and}\hspace{0.4cm}
		0 \to \uHom(P,M) \to \uHom(M_2,M) \to \uHom(M_1,M) \to 0
		\]
		are again exact. A diagram-chase using the nine-lemma shows that~$\dP$ is closed under extensions. Since~$\1\in\dP$ by assumption, we have~$\dT\subseteq\dP$.
		
		This implies in particular that every object of~$\dT$ is flat, and hence~$\dT$ is closed under tensor products, i.e.\ is a $\otimes$-subcategory of~$\dA$. Similarly, the functor~$\uHom(E,-)$ is exact for every~$E\in\dT$, which implies that~$\dT$ is also closed under~$\uHom(-,-)$, i.e.\ is a closed $\otimes$-subcategory.
		
		Finally, an object $E\in\dA$ is dualisable if and only if the natural maps
		\[
		\beta_{E,M}\colon \uHom(E,\1)\otimes M \to \uHom(E,M)
		\]
		are isomorphisms for all~$M\in\dA$. For fixed~$M$, the functors $\uHom(-,\1)\otimes M$ and~$\uHom(-,M)$ are exact when restricted to~$\dT$, and~$\beta_{\1,M}$ is an isomorphism for all~$M$. Hence by the five-lemma $\beta_{E,M}$ is an isomorphism for all~$E\in\dT$, i.e.\ all elements of~$\dT$ are dualisable.
		
		For~\eqref{condn:tannakian_lax_fibre_functor}, the maps
		\[
		\mu_{E_1,E_2}\colon \omega(E_1)\otimes\omega(E_2) \to \omega(E_1\otimes E_2)
		\]
		coming from the lax $\otimes$-structure on~$\omega$ are the components of a natural transformation between exact bifunctors. Since~$\mu_{\1,\1}$ is an isomorphism by assumption, a repeated application of the five-lemma shows that~$\mu_{E_1,E_2}$ is an isomorphism for all~$E_1,E_2$, i.e.\ $\omega$ is a $\otimes$-functor.
		
		To show that~$\omega$ is faithful for~$\dB\neq0$, consider the bifunctor $K\colon\dT^\op\times\dT \to \dB$ given by
		\[
		K(E_1,E_2) = \ker\left(\Hom_{\dT}(E_1,E_2) \to \Hom_F(\omega(E_1),\omega(E_2))\right) \,.
		\]
		It is easy to see that~$K$ is left-exact in each variable and $K(\1,\1)=0$ by~\eqref{condn:tannakian_endos_of_unit}, so by induction we have $K(E_1,E_2)=0$ for all~$E_1,E_2$ and we are done.
	\end{proof}
\end{lem}

\begin{ex}
	Suppose that~$(\dA,\otimes,\1)$ is the category of~$\cO_\fX$-modules on a ringed site~$(\fX,\cO_\fX)$ where~$\cO_\fX$ is a sheaf of~$F$-algebras. Then condition~\eqref{condn:tannakian_endos_of_unit} holds if and only if $\rH^0(\fX,\cO_\fX)=F$, and condition~\eqref{condn:tannakian_closed_tensor} holds automatically since every surjection~$M\twoheadrightarrow\cO_\fX$ is locally split. Given a morphism $x\colon(\ast,F)\to(\fX,\cO_\fX)$ of $F$-ringed sites, the pullback functor along~$x$ is an exact $\otimes$-functor which is faithful by condition~\eqref{condn:tannakian_lax_fibre_functor}. So, provided such an~$x$ exists and~$\fX$ is essentially small, the category of unipotent $\cO_\fX$-modules is a unipotent neutral Tannakian category.
\end{ex}

One advantage of finitely generated unipotent Tannakian categories~$\dT$ is that the universal object~${}_xE^\dT\in\pro{-}\dT$ pro-representing a fibre functor~$\omega_x$ admits an explicit description in terms of extensions. We define a sequence~$({}_xE^\dT_n)_{n\geq0}$ of objects of~$\dT$ as follows. We take~${}_xE^\dT_0=\1$, and for~$n\geq0$ define~${}_xE^\dT_{n+1}$ to be the extension
\begin{equation}\label{eq:tannakian_universal_extension}
	0 \to \1\otimes_F\Ext^1_\dT({}_xE^\dT_n,\1)^\vee \to {}_xE^\dT_{n+1} \to {}_xE^\dT_n \to 0
\end{equation}
corresponding to $1 \in \Ext^1_\dT({}_xE^\dT_n,\1\otimes_F\Ext^1_\dT({}_xE^\dT_n,\1)^\vee) = \End(\Ext^1_\dT({}_xE^\dT_n,\1))$. In other words,~\eqref{eq:tannakian_universal_extension} is the unique extension such that the coboundary map
\[
\delta\colon \Ext^1_\dT({}_xE^\dT_n,\1) \to \Ext^1_\dT({}_xE^\dT_n,\1)
\]
obtained by applying the functor~$\Hom_{\dT}(-,\1)$ to~\eqref{eq:tannakian_universal_extension} is the identity. We fix points~$e_{n,x}^\dT\in{}_xE^\dT_{n,x}$ by taking~$e_{0,x}^\dT=1$ and thereafter taking~$e_{n+1,x}^\dT$ to be any lift of~$e_{n,x}^\dT$ via the map in~\eqref{eq:tannakian_universal_extension}. This construction gives an explicit construction of~${}_xE^\dT$, as follows.

\begin{prop}\label{prop:tannakian_explicit_universal}
	Suppose that~$\dT$ is a finitely generated unipotent Tannakian category and that~$\omega_x$ is a fibre functor. Then for all~$n\geq0$, the pair~$({}_xE^\dT_n,e_{n,x}^\dT)$ represents the restriction of~$\omega_x$ to the category of objects of~$\dT$ of unipotency class~$\leq n$. In particular, $({}_xE^\dT,e_x^\dT)\colonequals (\varprojlim_n({}_xE^\dT_n),(e_{n,x}^\dT)_{n\geq0})$ pro-represents~$\omega_x$.
	\begin{proof}
		This is standard, see e.g.\ \cite[Theorem~2.1]{hadian:motivic_pi_1} or \cite[Proposition~3.3]{andreatta-iovita-kim}.
	\end{proof}
\end{prop}

\begin{rmk}
	The construction of the~${}_xE^\dT_n$ above shows that they, and~${}_xE^\dT$, are independent of the choice of fibre functor~$\omega_x$ up to non-canonical isomorphism. The reason we keep the~$x$ in the notation is that we will later want to endow certain universal objects~${}_xE^\dT$ with extra structures which do genuinely depend on the basepoint~$x$.
\end{rmk}

\subsection{A criterion for isomorphy}\label{ss:tannakian_comparison}

Suppose now that~$F'/F$ is an extension of fields, that~$\dT=(\dT,\otimes,\1)$ and~$\dT'=(\dT',\otimes',\1')$ are Tannakian categories over~$F$ and~$F'$, respectively, and that~$G\colon\dT\to\dT'$ is an $F$-linear exact $\otimes$-functor. We say that neutral fibre functors~$\omega_x$ and~$\omega_{x'}$ on~$\dT$ and~$\dT'$ respectively are \emph{compatible} just when the square
\begin{center}
\begin{tikzcd}[column sep = large]
	\dT \arrow[r,"G"]\arrow[d,"\omega_x"] & \dT' \arrow[d,"\omega_{x'}"] \\
	\Vec(F) \arrow[r,"F'\otimes_F(-)"] & \Vec(F')
\end{tikzcd}
\end{center}
commutes (up to a specified $\otimes$-natural isomorphism). If~$(\omega_x,\omega_{x'})$ is a compatible pair of neutral fibre functors, then there is an induced map
\begin{equation}\label{eq:tannakian_functoriality_family}
{}_{x'}E^{\dT'} \to G({}_xE^\dT)
\end{equation}
in~$\pro{-}\dT'$, namely the unique map sending~$e_{x'}^{\dT'}$ to the element~$1\otimes e_x^\dT\in G({}_xE^\dT)_{x'}\cong F'\otimes_F{}_xE^\dT_x$. And if~$(\omega_y,\omega_{y'})$ is a second compatible pair of fibre functors, then there is an induced morphism
\begin{equation}\label{eq:tannakian_functoriality}
\pi_1(\dT';\omega_{x'},\omega_{y'}) \to F'\otimes_F\pi_1(\dT;\omega_x,\omega_y)
\end{equation}
of affine~$F'$-schemes, given as follows. If~$\Lambda$ is an $F'$-algebra and~$\gamma\in\pi_1(\dT';\omega_{x'},\omega_{y'})(\Lambda)$, then~\eqref{eq:tannakian_functoriality} sends~$\gamma$ to the $\otimes$-natural transformation with components given by the composites.
\[
\Lambda\otimes_FE_x \cong \Lambda\otimes_{F'}G(E)_{x'} \xrightarrow{\gamma_{G(E)}} \Lambda\otimes_{F'}G(E)_{y'} \cong \Lambda\otimes_FE_y
\]
for~$E\in\dT$. We note that all of the obvious compatibilities hold.

\begin{lem}
	In the above setup, \eqref{eq:tannakian_functoriality_family} is compatible with coalgebra structures, and \eqref{eq:tannakian_functoriality} is compatible with the identities, composition and inversion in the fundamental groupoids. Moreover, the maps~\eqref{eq:tannakian_functoriality_family} and~\eqref{eq:tannakian_functoriality} are compatible with the identification in Proposition~\ref{prop:tannakian_groupoid_via_universal_objects}, in that the square
	\begin{center}
	\begin{tikzcd}
		\pi_1(\dT';\omega_{x'},\omega_{y'}) \arrow[r,"\eqref{eq:tannakian_functoriality}"]\arrow[d,"\eqref{eq:tannakian_groupoid_via_universal_objects}"',"\wr"] & F'\otimes_F\pi_1(\dT;\omega_x,\omega_y) \arrow[d,"\eqref{eq:tannakian_groupoid_via_universal_objects}"',"\wr"] \\
		\Spec({}_{x'}E^{\dT',\vee}_{y'}) \arrow[r,"\eqref{eq:tannakian_functoriality_family}"] & F'\otimes_F\Spec({}_xE^{\dT,\vee}_y)
	\end{tikzcd}
	\end{center}
	commutes.
	\begin{proof}
		Compatibility with coalgebra structures and groupoid operations is clear from the definitions. For the latter part, it suffices to show that the square
		\begin{equation}\label{eq:tannakian_compatibility_square}\tag{$\ast$}
		\begin{tikzcd}
			\Hom(\omega_{x'},\omega_{y'}) \arrow[r,"G_*"]\arrow[d,"\wr"] & F'\otimes_F\Hom(\omega_x,\omega_y) \arrow[d,"\wr"] \\
			{}_{x'}E^{\dT'}_{y'} \arrow[r,"\eqref{eq:tannakian_functoriality_family}"] & F'\otimes_F{}_xE^\dT_y
		\end{tikzcd}
		\end{equation}
		commutes, in which the vertical maps are the Yoneda isomorphisms~\eqref{eq:tannakian_yoneda} and the top map~$G_*$ is defined similarly to~\eqref{eq:tannakian_functoriality}. For any~$\gamma\in\Hom(\omega_{x'},\omega_{y'})$ we consider the diagram
		\begin{center}
		\begin{tikzcd}
			{}_{x'}E^{\dT'}_{x'} \arrow[r,"\eqref{eq:tannakian_functoriality_family}"]\arrow[d,"\gamma_{{}_{x'}E^{\dT'}}"] & G({}_xE^\dT)_{x'} \arrow[r,equals,"\sim"]\arrow[d,"\gamma_{G({}_xE^\dT)}"] & F'\otimes_F{}_xE^\dT_x \arrow[d,"G_*(\gamma)_{{}_xE^\dT}"] \\
			{}_{x'}E^{\dT'}_{y'} \arrow[r,"\eqref{eq:tannakian_functoriality_family}"] & G({}_xE^\dT)_{y'} \arrow[r,equals,"\sim"] & F'\otimes_F{}_xE^\dT_y \,,
		\end{tikzcd}
		\end{center}
		in which the left-hand square commutes by naturality and the right-hand square commutes by definition of~$G_*$. Pushing the element~$e_{x'}^{\dT'}$ both ways around this diagram shows that~\eqref{eq:tannakian_compatibility_square} commutes when evaluated at~$\gamma$, so we are done.
	\end{proof}
\end{lem}

We will need the following criterion for the maps~\eqref{eq:tannakian_functoriality_family} and~\eqref{eq:tannakian_functoriality} to be isomorphisms.

\begin{prop}\label{prop:tannakian_iso}
	Suppose that the Tannakian categories~$\dT$ and~$\dT'$ above are finitely generated unipotent, and that the map
	\[
	F'\otimes_F\Ext^1_\dT(\1,E) \to \Ext^1_{\dT'}(\1',G(E))
	\]
	induced by~$G$ is an isomorphism for all~$E\in\dT$. Then the maps~\eqref{eq:tannakian_functoriality_family} and~\eqref{eq:tannakian_functoriality} are isomorphisms for all compatible pairs of neutral fibre functors~$(\omega_x,\omega_{x'})$ and~$(\omega_y,\omega_{y'})$.
	\begin{proof}
		It suffices to prove that the maps
		\begin{equation}\label{eq:tannakian_functoriality_family_levelwise}
			\phi_n\colon{}_{x'}E^{\dT'}_n \to G({}_xE^\dT_n)
		\end{equation}
		sending~$e_{n,x'}^{\dT'}\in{}_{x'}E^{\dT'}_{n,x'}$ to~$1\otimes e_{n,x}^\dT\in G({}_xE^{\dT})_{x'}\cong F'\otimes_F{}_xE^\dT_x$ are isomorphisms for all~$n$. We do this using the explicit description of \S\ref{ss:tannakian_unipotent}, proceeding by induction. In the case~$n=0$ this is trivial, since~$\phi_0$ is the identity map on~$\1'$.
		
		For~$n\geq0$, suppose that~$\phi_n$ is an isomorphism. It follows by definition that~$\phi_{n+1}$ fits into a commuting diagram
		\begin{equation}\label{diag:comparison_of_exact_seqs}
		\begin{tikzcd}
			0 \arrow[r] & \1'\otimes_{F'}\Ext^1_{\dT'}({}_{x'}E^{\dT'}_n,\1')^\vee \arrow[r]\arrow[d] & {}_{x'}E^{\dT'}_{n+1} \arrow[r]\arrow[d,"\phi_{n+1}"] & {}_{x'}E^{\dT'}_n \arrow[r]\arrow[d,"\phi_n","\wr"'] & 0 \\
			0 \arrow[r] & \1'\otimes_F\Ext^1_{\dT}({}_xE^\dT_n,\1)^\vee \arrow[r] & G({}_xE^\dT_{n+1}) \arrow[r] & G({}_xE^\dT_n) \arrow[r] & 0
		\end{tikzcd}
		\end{equation}
		where the rows are the exact sequences~\eqref{eq:tannakian_universal_extension}. Applying the functor~$\Hom_{\dT'}(-,\1')$, we obtain a commuting square
		\begin{equation}\label{diag:tannakian_coboundary_square}
		\begin{tikzcd}
			\Ext^1_{\dT'}({}_{x'}E^{\dT'}_n,\1') \arrow[r,"\delta'"] & \Ext^1_{\dT'}({}_{x'}E^{\dT'}_n,\1') \\
			F'\otimes_F\Ext^1_{\dT}({}_xE^\dT_n,\1) \arrow[r,"\delta"]\arrow[u] & \Ext^1_{\dT'}(G({}_xE^\dT_n),\1') \arrow[u,"\phi_n^*"',"\wr"] \,.
		\end{tikzcd}
		\end{equation}
		The topmost coboundary map~$\delta'$ is the identity map by the description of the coboundary map in~\eqref{eq:tannakian_universal_extension}, while the bottom coboundary maps~$\delta$ is the map induced by~$G$, which is an isomorphism by assumption. Hence the left-hand map in~\eqref{diag:tannakian_coboundary_square} is also an isomorphism. This implies that the left-hand vertical map in~\eqref{diag:comparison_of_exact_seqs} is also an isomorphism, and hence~$\phi_{n+1}$ is an isomorphism by the five-lemma. This completes the induction.
	\end{proof}
\end{prop}

\begin{rmk}
	In the particular case that~$F'=F$, Proposition~\ref{prop:tannakian_iso} recovers the well-known criterion that an exact $F$-linear $\otimes$-functor between unipotent neutral Tannakian categories over~$F$ is an equivalence if and only if it induces an isomorphism on first Ext-groups.
\end{rmk}

\section{The unipotent Riemann--Hilbert fundamental groupoid}

Now we come to the central construction in this paper: the \emph{pro-unipotent Riemann--Hilbert fundamental groupoid} of a smooth geometrically connected variety~$Y/K$, which is a pro-unipotent groupoid-scheme over~$\sB_\deR$ which we will use to relate the pro-unipotent \'etale and de Rham fundamental groupoids of~$Y$. The construction follows the methods of $p$-adic Hodge theory for log rigid analytic varieties and the $p$-adic Riemann--Hilbert correspondence of~\cite{dllz:riemann-hilbert}.

As in the introduction, let~$Y/K$ be a smooth geometrically connected variety, and choose a smooth proper compactification~$X$ of~$Y$ whose complement~$D$ is a normal crossings divisor. Let~$X^\an$ denote the rigid analytification of~$X$, viewed as an adic space in the usual way, and make~$X^\an$ into an fs log adic space by giving it the log structure coming from the divisor~$D^\an$, as in \cite[Example~2.1.2]{dllz:riemann-hilbert}.

Let~$\cO_{X^\an}$ denote the structure sheaf on the analytic site~$X^\an_\an$ of~$X^\an$, and let~$\cO_\dX\colonequals\cO_{X^\an}\hatotimes_K\sB_\deR$ be the completion of~$\cO_{X^\an}\otimes_K\sB_\deR$ with respect to its filtration, i.e.
\[
\cO_\dX = \left(\varprojlim_n \cO_{X^\an}\otimes_K(\sB_\deR^+/t^n\sB_\deR^+)\right)[t^{-1}] \,.
\]
This is simultaneously a sheaf of~$\sB_\deR$-algebras and a sheaf of $\cO_{X^\an}$-algebras on~$X^\an_\an$, and comes with a $\sB_\deR$-linear logarithmic connection
\[
\nabla\colon \cO_\dX \to \cO_\dX \otimes_{\cO_{X^\an}}\Omega^1_{X^\an/K}
\]
compatible with the derivation on~$\cO_{X^\an}$, where~$\Omega^1_{X^\an/K}$ is the sheaf of logarithmic differentials on~$X^\an$.

We write~$\dX$ for the ringed site~$(X^\an_\an,\cO_\dX)$. By a \emph{connection} on an~$\cO_\dX$-module~$\dE$, we mean a morphism
\[
\nabla\colon \dE \to \dE\otimes_{\cO_{X^\an}}\Omega^1_{X^\an/K}
\]
of abelian sheaves satisfying the Leibniz rule with respect to the connection on~$\cO_\dX$. In particular, such a connection is always $\sB_\deR$-linear. We say that~$\nabla$ is \emph{integrable} just when it is integrable in the sense of connections on $\cO_{X^\an}$-modules.

We write~$\MIC^\un(\dX,\cO_\dX)$ for the category of unipotent objects in the category of $\cO_\dX$-modules with integrable connection. Objects of~$\MIC^\un(\dX,\cO_\dX)$ are called \emph{unipotent $\cO_\dX$-vector bundles with integrable connection}; as the name suggests, they are automatically locally free as $\cO_\dX$-modules. Given a point~$x\in Y^\an(L)$ defined over a finite extension~$L$ of~$K$ inside~$\bC_K$, the natural morphism $(\ast,L)\to(X^\an_\an,\cO_{X^\an})$ of ringed sites extends naturally to a morphism
\[
i_x\colon (\ast,\sB_\deR) \to (\dX,\cO_\dX) \,.
\]
We write
\[
\omega_x^\RH\colon \MIC^\un(\dX,\cO_\dX) \to \Mod(\sB_\deR)
\]
for the map induced by the pullback functor~$i_x^*$.

\begin{prop}\label{prop:RH_is_tannakian}
	$\MIC^\un(\dX,\cO_\dX)$ is a finitely generated unipotent Tannakian category over~$\sB_\deR$ (with respect to the evident tensor product), and~$\omega^\RH_x$ is a neutral fibre functor for all points~$x\in Y^\an(L)$ as above.
\end{prop}

The subtle part of this proposition is the finite generation. For this, we need to have a good description of the Ext-groups in~$\MIC^\un(\dX,\cO_\dX)$. This is afforded by de Rham cohomology, as usual.

\begin{lem}\label{lem:RH_extensions}
	For any~$\dE_1,\dE_2\in\MIC^\un(\dX,\cO_\dX)$ there is a canonical $\sB_\deR$-linear isomorphism
	\[
	\Ext^1_{\MIC}(\dE_1,\dE_2) \cong \rH^1_\deR(\dX,\dE_1^\vee\otimes_{\cO_\dX}\dE_2) \,,
	\]
	natural in~$\dE_1$ and~$\dE_2$. Here $\rH^1_\deR(\dX,\dE)\colonequals\bH^1(X^\an_\an,\DeR(\dE))$, where~$\DeR(\dE)$ is the de Rham complex of~$\dE$.
	\begin{proof}
		This is a variation on the usual argument, see e.g.\ \cite[Proposition~3.3]{hadian:motivic_pi_1}. Let~$\dE$ be an extension of~$\dE_1$ by~$\dE_2$. Choose a covering~$\fU=(U_i)_{i\in I}$ of~$\dX=X^\an_\an$ over which the extension splits $\cO_\dX$-linearly, and fix a choice of $\cO_\dX$-linear splitting~$s_i$ of~$\dE|_{U_i}\to\dE_1|_{U_i}$ for each~$i$. Viewing~$s_i$ as an element of~$\rH^0(U_i,\dE_1^\vee\otimes_{\cO_\dX}\dE)$, we obtain elements
		\[
		\kappa_{ij}\colonequals s_i|_{U_{ij}} - s_j|_{U_{ij}} \in \rH^0(U_{ij},\dE_1^\vee\otimes_{\cO_\dX}\dE_2)
		\hspace{0.4cm}\text{and}\hspace{0.4cm}
		\lambda_i\colonequals \nabla(s_i) \in \rH^0(U_i,\dE_1^\vee\otimes_{\cO_\dX}\dE_2\otimes_{\cO_{X^\an}}\Omega^1_{X^\an/K})
		\]
		for each~$i,j$. Taken together, these elements define a \u Cech $1$-cochain
		\[
		\xi=((\kappa_{ij})_{i,j\in I},(\lambda_i)_{i\in I}) \in \breve \rC^1\!(\fU,\DeR(\dE_1^\vee\otimes_{\cO_\dX}\dE_2)) \,.
		\]
		Since the connection~$\nabla$ on~$\dE$ is integrable, we have~$\rd\xi=0$, i.e.\ $\xi$ is a $1$-cocycle, and one easily checks that different choices of the splittings~$s_i$ only change~$\xi$ by a coboundary. Hence we have constructed a map
		\begin{equation}\label{eq:cohomology_classes_over_coverings}\tag{$\ast$}
			\{\text{extensions~$\dE$ that are $\cO_\dX$-linearly split over~$\fU$}\}/\text{iso} \to \breve \bH^1\!(\fU,\DeR(\dE_1^\vee\otimes_{\cO_\dX}\dE_2))
		\end{equation}
		given by~$[\dE]\mapsto[\xi]$.
		
		To show that~\eqref{eq:cohomology_classes_over_coverings} is bijective, we explicitly describe the inverse. Given a \u Cech $1$-cocycle
		\[
		\xi=((\kappa_{ij})_{i,j\in I},(\lambda_i)_{i\in I}) \in \breve \rC^1\!(\fU,\DeR(\dE_1^\vee\otimes_{\cO_\dX}\dE_2)) \,,
		\]
		we consider the connections~$\nabla_i$ on the $\cO_\dX|_{U_i}$-modules $\dE_2|_{U_i}\oplus\dE_1|_{U_i}$ given by
		\[
		\nabla_i \colonequals 
		\begin{pmatrix}
			\nabla_{\dE_2}|_{U_i} & \lambda_i \\
			0 & \nabla_{\dE_1}|_{U_i}
		\end{pmatrix} \,.
		\]
		Since~$\xi$ is a $1$-cocycle, we have~$\nabla(\lambda_i)=0$ for each~$i$, so~$\nabla_i$ is an integrable connection. Moreover, we have~$\nabla(\kappa_{ij})=\lambda_i-\lambda_j$ for all~$i,j$, which implies the automorphisms
		\[
		\phi_{ij}\colonequals
		\begin{pmatrix}
			1_{\dE_2}|_{U_{ij}} & \kappa_{ij} \\
			0 & 1_{\dE_1}|_{U_{ij}}
		\end{pmatrix}
		\colon \dE_2|_{U_{ij}}\oplus\dE_1|_{U_{ij}} \xrightarrow\sim \dE_2|_{U_{ij}}\oplus\dE_1|_{U_{ij}}
		\]
		carry~$\nabla_i|_{U_{ij}}$ to~$\nabla_j|_{U_{ij}}$. Hence if~$\dE$ is the $\cO_\dX$-linear extension of~$\dE_1$ by~$\dE_2$ obtained by gluing the~$\dE_2|_{U_i}\oplus\dE_1|_{U_i}$ along the isomorphisms~$\phi_{ij}$, then the connections~$\nabla_i$ induce a well-defined connection on~$\dE$, making it into an extension of~$\dE_1$ by~$\dE_2$ in~$\MIC^\un(\dX\,\cO_\dX)$.
		
		This construction is easily checked to define an inverse to~\eqref{eq:cohomology_classes_over_coverings}. Taking the colimit over coverings~$\fU$ then yields the desired result.
	\end{proof}
\end{lem}

We have the following calculation of these cohomology groups in the case~$\dE=\cO_\dX$.

\begin{lem}\label{lem:RH-dR_extension_comparison_baby}
	The natural map
	\[
	\sB_\deR\otimes_K\rH^1_\deR(X^\an,\cO_{X^\an}) \to \rH^1_\deR(\dX,\cO_\dX)
	\]
	is an $\sB_\deR$-linear isomorphism. In particular, $\rH^1_\deR(\dX,\cO_\dX)$ is finite-dimensional over~$\sB_\deR$.
	\begin{proof}
		See~\cite[Lemma~3.6.3]{dllz:riemann-hilbert}. We note for later reference that the proof of \cite[Lemma~3.6.3]{dllz:riemann-hilbert} establishes a stronger result, that the natural map
		\[
		\sB_\deR\otimes_K\rR\!\Gamma(X^\an_\an,\DeR(\cO_{X^\an})) \to \rR\!\Gamma(\dX,\DeR(\cO_\dX))
		\]
		is a filtered quasi-isomorphism, i.e.\ is an isomorphism in the filtered derived category.
	\end{proof}
\end{lem}


\begin{proof}[Proof of Proposition~\ref{prop:RH_is_tannakian}]
	The criteria of Lemma~\ref{lem:tannakian_criteria} are easily verified (with~$\dA$ the $\sB_\deR$-linear abelian closed $\otimes$-category of all $\cO_\dX$-modules with integrable connection). Finite-dimensionality of $\Ext^1_{\MIC}(\cO_\dX,\cO_\dX)$ follows from Lemmas~\ref{lem:RH_extensions} and~\ref{lem:RH-dR_extension_comparison_baby}.
\end{proof}

As a consequence, we can now make the following definition.

\begin{defi}
	For two points~$x,y\in Y^\an(L)$ defined over a finite extension~$L$ of~$K$ inside~$\bC_K$, we define the \emph{pro-unipotent Riemann--Hilbert path torsor}
	\[
	\pi_1^\RH(\dX;x,y)
	\]
	to be the Tannakian path torsor from~$\omega_x^\RH$ to~$\omega_y^\RH$ in the Tannakian category $\MIC^\un(\dX,\cO_\dX)$. These assemble into a groupoid in affine $\sB_\deR$-schemes on object-set~$Y(\Kbar)$, which we call the \emph{pro-unipotent Riemann--Hilbert fundamental groupoid of~$X^\an_{\bC_K}$}.
\end{defi}

This fundamental groupoid comes with extra structures, which we now describe.

\subsection{Galois action}

The left action of~$G_K$ on~$\sB_\deR$ induces an action on~$\cO_\dX$, and hence a right action on the category~$\MIC^\un(\dX,\cO_\dX)$. Explicitly, if~$\dE\in\MIC^\un(\dX,\cO_\dX)$ and $\sigma\in G_K$, then~$\sigma^*(\dE)$ has the same underlying abelian sheaf and connection as~$\dE$, but has the $\sigma$-twisted $\cO_\dX$-action:
\[
\lambda\colon e \mapsto \sigma(\lambda)e
\]
for sections~$\lambda$ of~$\cO_\dX$ and~$e$ of~$\dE$.

If~$x\in Y^\an(L)$ is an $L$-rational point, then its associated fibre functor~$\omega_x^\RH\colon\MIC^\un(\dX,\cO_\dX)\to\Vec(\sB_\deR)$ is $G_L$-equivariant. Hence for any second $L$-rational point~$y\in Y^\an(L)$, the functor of $\otimes$-natural isomorphisms $\omega_x^\RH\xrightarrow\sim\omega_y^\RH$ is $G_L$-invariant, and so there is an induced left $G_L$-action on~$\pi_1^\RH(\dX;x,y)$. This action is semilinear, meaning that its action on the affine ring of~$\pi_1^\RH(\dX;x,y)$ is semilinear. It is clear from this construction that the $G_L$-actions on the schemes $\pi_1^\RH(\dX;x,y)$ for varying~$x,y\in Y^\an(L)$ respect the groupoid structure maps.


Similarly, there is a semilinear action of~$G_L$ on the universal object~${}_x\dE^\RH$, where
\[
\phi_\sigma\colon {}_x\dE^\RH \to \sigma^*({}_x\dE^\RH) = {}_x\dE^\RH
\]
is the unique map sending~$e_x^\RH$ to $e_x^\RH\in(\sigma^*({}_x\dE^\RH))_x\cong\sigma^*({}_x\dE^\RH_x)={}_x\dE^\RH_x$. This action can be characterised by a universal property.

\begin{lem}\label{lem:RH_equivariant_representability}
	Let~$\MIC^\un(\dX,\cO_\dX)^{G_L}$ denote the category of unipotent $\cO_\dX$-vector bundles with integrable connection endowed with a semilinear action of~$G_L$ preserving the connection. Then $({}_x\dE^\RH,e_x^\RH)$, endowed with the action described above, pro-represents the functor
	\begin{align*}
		\MIC^\un(\dX,\cO_\dX)^{G_L} &\to \Mod(\sB_\deR^{G_L}) \\
		\dE &\mapsto (\dE_x)^{G_L}
	\end{align*}
	\begin{proof}
		Since $({}_x\dE^\RH,e^\RH_x)$ pro-represents the fibre functor
		\[
		\omega^\RH_x\colon\MIC^\un(\dX,\cO_\dX)\to\Vec(\sB_\deR) \,,
		\]
		the proposition amounts to the following: for any~$\dE\in\MIC^\un(\dX,\cO_\dX)^{G_L}$, a morphism $f\colon{}_x\dE^\RH\to\dE$ in $\pro{-}\MIC^\un(\dX,\cO_\dX)$ is $G_L$-equivariant if and only if~$f(e^\RH_x)\in(\dE_x)^{G_L}$. But this is straightforward: $f$ is $G_L$-equivariant if and only if the square
		\begin{center}
			\begin{tikzcd}
				{}_x\dE^\RH \arrow[r,"\phi_\sigma"]\arrow[d,"f"] & \sigma^*({}_x\dE^\RH) \arrow[d,"\sigma^*(f)"] \\
				\dE \arrow[r,"\phi_\sigma"] & \sigma^*(\dE)
			\end{tikzcd}
		\end{center}
		commutes for all~$\sigma\in G_L$ (where~$\phi_\sigma$ are the actions of~$\sigma$ on~${}_x\dE^\RH$ and~$\dE$). Each such square commutes if and only if it commutes when evaluated at~$e_x^\RH\in{}_x\dE^\RH_x$, whence the result.
	\end{proof}
\end{lem}

Let us also remark that the $G_L$-actions on the universal objects~${}_x\dE^\RH$ and the Tannakian path-torsors $\pi_1^\RH(\dX;x,y)$ are compatible, in the sense that the latter is obtained from taking the fibre of the former at~$y$ and applying Proposition~\ref{prop:tannakian_groupoid_via_universal_objects}. This is a consequence of the following lemma.

\begin{lem}\label{lem:RH_yoneda_equivariant}
	For all~$x,y\in Y^\an(L)$, the Yoneda isomorphism~\eqref{eq:tannakian_yoneda}
	\[
	\Hom(\omega^\RH_x,\omega^\RH_y) \cong {}_x\dE^\RH_y
	\]
	is $G_L$-equivariant, for the natural~$G_L$-action on either side.
	\begin{proof}
		For~$\gamma\in\Hom(\omega^\RH_x,\omega^\RH_y)$, $\sigma(\gamma)$ sends~$e^\RH_x\in{}_x\dE^\RH_x$ to~$\sigma(\gamma_{\sigma^*({}_x\dE^\RH)}(\phi_\sigma(e^\RH_x)))=\sigma(\gamma_{{}_x\dE^\RH}(e^\RH_x))$, since~$e^\RH_x$ is~$G_L$-invariant. The result follows.
	\end{proof}
\end{lem}

\subsection{Hodge filtration}\label{ss:RH_hodge_filtration}

The other structure on the Riemann--Hilbert fundamental groupoid is a Hodge filtration. The usual $t$-adic filtration on~$\sB_\deR$ induces a filtration on~$\cO_\dX$ which satisfies Griffiths transversality with respect to the connection.

\begin{defi}\label{def:RH_filtrations}
	Let~$(\dE,\nabla)$ be an $\cO_\dX$-module with integrable connection. By a \emph{filtration} on~$\dE$ we mean a decreasing filtration
	\[
	\dots\geq\rF^1\!\dE\geq\rF^0\!\dE\geq\rF^{-1}\!\dE\geq\dots
	\]
	of~$\dE$ which is compatible with the filtration on~$\cO_\dX$ and satisfies Griffiths transversality with respect to the connection~$\nabla$. We define the category
	\[
	\FMIC^\un(\dX,\cO_\dX)
	\]
	of \emph{filtered $\cO_\dX$-vector bundles with integrable connection} to be the category of unipotent objects in the category of filtered $\cO_\dX$-modules with integrable connection. The filtration on an object of~$\FMIC^\un(\dX,\cO_\dX)$ is automatically exhaustive and separated, and each~$\rF^i\!\dE$ is always a local direct summand in~$\dE$.
\end{defi}

We will continually make use of the following observation: if~$\rF^\bullet$ is a filtration on a~$\cO_\dX$-module~$\dE$ which is compatible with the filtration on~$\cO_\dX$, then~$\rF^i\!\dE=t^i\rF^0\!\dE$ for all~$i$. So a filtration is determined by~$\rF^0$.


Now fix a point~$x\in Y^\an(L)$ defined over a finite extension of~$K$, and let $(\dE^\RH,e^\RH)=({}_x\dE^\RH,e_x^\RH)$ denote the pair pro-representing the fibre functor~$\omega^\RH_x$. We will construct a Hodge filtration on~$\dE^\RH$, making it into a pro-object of~$\FMIC^\un(\dX,\cO_\dX)$. This is done using the recursive description of~$\dE^\RH$ from Proposition~\ref{prop:tannakian_explicit_universal} as follows.

\begin{prop}\label{prop:RH_hodge_filtration}
	There is a unique way to put a filtration~$\rF^\bullet$ on each~$\dE^\RH_n$ such that the following conditions are satisfied:
	\begin{enumerate}\setcounter{enumi}{-1}
		\item\label{condn:RH_hodge_basic} $\dE^\RH_n\in\FMIC^\un(\dX,\cO_\dX)$ for all~$n$;
		\item\label{condn:RH_hodge_normalisation} the filtration on~$\dE^\RH_0=\cO_\dX$ is the $t$-adic filtration;
		\item\label{condn:RH_hodge_exactness} for all~$n\geq1$ the sequence
		\begin{equation}\label{eq:RH_universal_extension}
			0 \to \cO_\dX\otimes_{\sB_\deR}\rH^1_\deR(\dX,\dE^{\RH,\vee}_{n-1})^\vee \to \dE^\RH_n \to \dE^\RH_{n-1} \to 0
		\end{equation}
		is strict exact for the filtrations; and
		\item\label{condn:RH_hodge_rigidity} $e^\RH_n\in\rF^0\!\dE^\RH_{n,x}$ for all~$n$.
	\end{enumerate}
	The sequence~\eqref{eq:RH_universal_extension} in the second point is the sequence~\eqref{eq:tannakian_universal_extension}, using the description of the Ext-groups from Lemma~\ref{lem:RH_extensions}. The implicit filtration on~$\rH^1_\deR(\dX,\dE^{\RH,\vee}_{n-1})$ is the one induced from the filtration on the de Rham complex induced from the filtration on~$\dE^\RH_{n-1}$.
\end{prop}

In preparation for this result, we first need an explicit description of the groups of strict extensions in~$\FMIC^\un(\dX,\cO_\dX)$, which we denote by~$\Ext^1_{\FMIC}(-,-)$. We define for any~$\dE\in\FMIC^\un(\dX,\cO_\dX)$
\[
\rH^1_\FdeR(\dX,\dE) \colonequals \bH^1(X^\an_\an,\rF^0\!\DeR(\dE)) \,,
\]
where~$\rF^\bullet\!\DeR(\dE)$ denotes the induced filtration on the de Rham complex of~$\dE$.

\begin{lem}\label{lem:filtered_RH_extensions}
	For any~$\dE_1,\dE_2\in\FMIC^\un(\dX,\cO_\dX)$ there is a canonical $\sB_\deR^+$-linear isomorphism
	\[
	\Ext^1_{\FMIC}(\dE_1,\dE_2) \cong \rH^1_{\FdeR}(\dX,\dE_1^\vee\otimes_{\cO_\dX}\dE_2) \,,
	\]
	natural in~$\dE_1$ and~$\dE_2$.
	\begin{proof}
		A very similar proof to that of Lemma~\ref{lem:RH_extensions} works here. Given a strict extension~$\dE$ of~$\dE_1$ by~$\dE_2$, one may choose a covering~$\fU=(U_i)_{i\in I}$ over which~$\dE$ splits as an extension of filtered $\cO_\dX$-modules. If one chooses splittings~$s_i$ of~$\dE|_{U_i}\to\dE_1|_{U_i}$ compatible with filtrations, then the elements $\kappa_{ij}\colonequals s_i|_{U_{ij}}-s_j|_{U_{ij}}$ and $\lambda_i\colonequals\nabla(s_i)$ lie in $\rH^0(U_{ij},\rF^0\!(\dE_1^\vee\otimes_{\cO_\dX}\dE_2))$ and $\rH^0(U_i,\rF^{-1}\!(\dE_1^\vee\otimes_{\cO_\dX}\dE_2)\otimes_{\cO_{X^\an}}\Omega^1_{X^\an/K})$, respectively, the latter by Griffiths transversality. So $\xi\colonequals((\kappa_{ij})_{i,j},(\lambda_i)_i)$ defines a \u Cech $1$-cochain with coefficients in~$\rF^0\!\DeR(\dE_1^\vee\otimes_{\cO_\dX}\dE_2)$. As in the proof of Lemma~\ref{lem:RH_extensions}, one checks that~$\xi$ is a $1$-cocycle, that it is independent of the choice of filtered splittings~$s_i$ up to $1$-coboundaries, and that the assignment~$[\dE]\mapsto[\xi]$ gives the desired bijection.
	\end{proof}
\end{lem}

\begin{lem}\label{lem:RH_filtered_extensions_finite_type}
	$\rH^1_\FdeR(\dX,\dE)$ is a $\sB_\deR^+$-module of finite type for all~$\dE\in\FMIC^\un(\dX,\cO_\dX)$.
	\begin{proof}
		It suffices by an inductive argument to establish it in the case~$\dE=\cO_\dX$. For this, the fact that the map
		\[
		\sB_\deR\otimes_K\rR\!\Gamma(X^\an_\an,\cO_{X^\an}) \to \rR\!\Gamma(\dX,\cO_\dX)
		\]
		is an isomorphism in the filtered derived category (see Lemma~\ref{lem:RH-dR_extension_comparison_baby}) implies that it induces an isomorphism on the cohomology of the $0$th filtered part. Concretely, this means that~$\rH^1_\FdeR(\dX,\cO_\dX)$ sits in a pushout square
		\begin{center}
			\begin{tikzcd}
				\sB_\deR^+\otimes_K\rH^1(X^\an_\an,\rF^1\!\DeR(\cO_{X^\an})) \arrow[r,hook]\arrow[d,hook] & \sB_\deR^+\otimes_K\rH^1(X^\an_\an,\rF^0\!\DeR(\cO_{X^\an})) \arrow[d] \\
				t^{-1}\sB_\deR^+\otimes_K\rH^1(X^\an_\an,\rF^1\!\DeR(\cO_{X^\an})) \arrow[r] & \rH^1_\FdeR(\dX,\cO_\dX) \,,
			\end{tikzcd}
		\end{center}
		where the injectivity of the top map follows from the degeneration of the Hodge--de Rham spectral sequence~$\rH^\bullet_\deR(X^\an,\cO_{X^\an})$ at the first page. The result follows.
	\end{proof}
\end{lem}

For any~$\dE\in\FMIC^\un(\dX,\cO_\dX)$, there is a natural map
\begin{equation}\label{eq:RH_forget_filtration}
\rH^1_\FdeR(\dX,\dE) \to \rH^1_\deR(\dX,\dE)
\end{equation}
induced by the inclusion~$\rF^0\!\DeR(\dE)\hookrightarrow\DeR(\dE)$; via the identifications in Lemmas~\ref{lem:RH_extensions} and~\ref{lem:filtered_RH_extensions} this corresponds to the map
\[
\Ext^1_{\FMIC}(\cO_\dX,\dE) \to \Ext^1_{\MIC}(\cO_\dX,\dE)
\]
given by forgetting the filtration. The map~\eqref{eq:RH_forget_filtration} has image~$\rF^0\!\rH^1_\deR(\dX,\dE)$ by definition, and induces a filtered isomorphism
\[
\rH^1_\FdeR(\dX,\dE)[t^{-1}] = \sB_\deR\otimes_{\sB_\deR^+}\rH^1_\FdeR(\dX,\dE) \xrightarrow\sim \rH^1_\deR(\dX,\dE) \,;
\]
in particular, $\rH^1_\deR(\dX,\dE)$ is filtered-isomorphic to a finite direct sum of copies of~$\sB_\deR$. Note however that the map~\eqref{eq:RH_forget_filtration} is not \emph{a priori} injective in general.

\begin{lem}\label{lem:OB_deR_good}
	\eqref{eq:RH_forget_filtration} is injective for~$\dE=\cO_\dX$.
	\begin{proof}
		The proof of Lemma~\ref{lem:RH_filtered_extensions_finite_type} showed that~$\rH^1_\FdeR(\dX,\cO_\dX)=\rF^0(\sB_\deR\otimes_K\rH^1_\deR(X^\an,\cO_{X^\an}))$, which injects into~$\rH^1_\deR(\dX,\cO_\dX)=\sB_\deR\otimes_K\rH^1_\deR(X^\an,\cO_{X^\an})$.
	\end{proof}
\end{lem}

%
%


\begin{proof}[Proof of Proposition~\ref{prop:RH_hodge_filtration}]
	We use a modification of the argument sketched in \cite[Lemma~3.6]{hadian:motivic_pi_1}. We will show by induction that there exists a unique sequence of filtrations on the~$\dE^\RH_n$ satisfying conditions \eqref{condn:RH_hodge_basic}--\eqref{condn:RH_hodge_rigidity}, and that these filtrations additionally satisfy the auxiliary condition
	\textit{
		\begin{enumerate}\setcounter{enumi}{3}
			\item\label{condn:RH_hodge_good} the map~\eqref{eq:RH_forget_filtration} $\rH^1_\FdeR(\dX,\dE^{\RH,\vee}_n)\to\rH^1_\deR(\dX,\dE^{\RH,\vee}_n)$ is injective.
		\end{enumerate}
	}
	In the base case~$n=0$, the given filtration on~$\dE^\RH_0=\cO_\dX$ is the unique filtration satisfying~\eqref{condn:RH_hodge_normalisation}. It clearly satisfies~\eqref{condn:RH_hodge_basic} and~\eqref{condn:RH_hodge_rigidity}, and satisfies~\eqref{condn:RH_hodge_good} by Lemma~\ref{lem:OB_deR_good}.
	
	So suppose inductively that there is a unique sequence of filtrations on~$\dE^\RH_m$ for~$0\leq m\leq n$ satisfying~\eqref{condn:RH_hodge_basic}--\eqref{condn:RH_hodge_rigidity} for~$m\leq n$, and that these filtrations also satisfy~\eqref{condn:RH_hodge_good}. We write $\dG:=\dE^{\RH,\vee}_n\otimes_{\sB_\deR}\rH^1_\deR(\dX,\dE^{\RH,\vee}_n)^\vee$, and endow it with the obvious filtration. Since $\dG$ is isomorphic to a direct sum of copies of $\dE^{\RH,\vee}_n$, it follows from~\eqref{condn:RH_hodge_good} that the map
	\[
	\rH^1_\FdeR(\dX,\dG)\to\rH^1_\deR(\dX,\dG)\cong\End(\rH^1_\deR(\dX,\dE^{\RH,\vee}_n))
	\]
	is injective. Since the class of the extension~$\dE^\RH_{n+1}$ is $1\in\rF^0\End(\rH^1_\deR(\dX,\dE^{\RH,\vee}_n))$, we see from Lemma~\ref{lem:filtered_RH_extensions} that there is a strict extension of~$\dE^\RH_n$ by $\cO_\dX\otimes_{\sB_\deR}\rH^1_\deR(\dX,\dE^{\RH,\vee}_n)^\vee$ whose underlying non-filtered extension class is~$[\dE^\RH_{n+1}]$, and this extension is unique up to isomorphism of strict extensions.
	
	In other words, there exists a filtration~$\rF^\bullet$ on~$\dE^\RH_{n+1}$ making~\eqref{eq:RH_universal_extension} into a strict extension, and this filtration is unique up to the action of (non-filtered) automorphisms of the extension~\eqref{eq:RH_universal_extension}. By Proposition~\ref{prop:tannakian_explicit_universal}, this automorphism group is
	\[
	\Hom_{\MIC}(\dE^\RH_n,\cO_\dX\otimes_{\sB_\deR}\rH^1_\deR(\dX,\dE^{\RH,\vee}_n)^\vee)=\rH^1_\deR(\dX,\dE^{\RH,\vee}_n)^\vee \,.
	\]
	The action of this automorphism group on the fibre of~$\dE^\RH_{n+1}$ at~$x$ is the addition action of the subgroup $\rH^1_\deR(\dX,\dE^{\RH,\vee}_n)^\vee\leq\dE^\RH_{n+1,x}$; in particular it acts simply transitively on the fibres of $\dE^\RH_{n+1,x}\to\dE^\RH_{n,x}$. Hence, by using a suitable automorphism, we may assume that~$e^\RH_{n+1}\in\rF^0\!\dE^\RH_{n+1,x}$.
	
	Thus the filtration on~$\dE^\RH_{n+1}$ constructed above satisfies conditions~\eqref{condn:RH_hodge_basic}--\eqref{condn:RH_hodge_rigidity}. To see that it is the unique such filtration, we observe that any other filtration satisfying~\eqref{condn:RH_hodge_basic}--\eqref{condn:RH_hodge_rigidity} would have to differ from~$\rF^\bullet$ by the action of an element $\alpha\in\rF^0\rH^1_\deR(\dX,\dE^{\RH,\vee}_{n,x})^\vee=\rF^0\!\Hom_{\MIC}(\dE^\RH_n,\cO_\dX\otimes_{\sB_\deR}\rH^1_\deR(\dX,\dE^{\RH,\vee}_{n,x})^\vee)$. So the action of~$\alpha$ would fix the filtration~$\rF^\bullet$, which is therefore unique.
	
	To complete the induction, it remains to show that this filtration on~$\dE^\RH_{n+1}$ also satisfies condition~\eqref{condn:RH_hodge_good}. For this, we take logarithmic de Rham cohomology of the dual of the strict exact sequence~\eqref{eq:RH_universal_extension} to obtain a commuting diagram
	{\small
		\begin{center}
			\begin{tikzcd}[column sep = small]
				\rF^0\rH^1_\deR(\dX,\dE^{\RH,\vee}_n) \arrow[r,"\delta","\sim"']\arrow[d,hook] & \rH^1_\FdeR(\dX,\dE^{\RH,\vee}_n) \arrow[r]\arrow[d,hook] & \rH^1_\FdeR(\dX,\dE^{\RH,\vee}_{n+1}) \arrow[r]\arrow[d] & \rF^0\!\left(\rH^1_\deR(\dX,\cO_\dX)\otimes\rH^1_\deR(\dX,\dE^{\RH,\vee}_n)\right) \arrow[d,hook] \\
				\rH^1_\deR(\dX,\dE^{\RH,\vee}_n) \arrow[r,"\delta","\sim"'] & \rH^1_\deR(\dX,\dE^{\RH,\vee}_n) \arrow[r] & \rH^1_\deR(\dX,\dE^{\RH,\vee}_{n+1}) \arrow[r] & \rH^1_\deR(\dX,\cO_\dX)\otimes\rH^1_\deR(\dX,\dE^{\RH,\vee}_n)
			\end{tikzcd}
		\end{center}
	}
	\noindent with exact rows. The final vertical map is injective since $\cO_\dX\otimes\rH^1_\deR(\dX,\dE^{\RH,\vee}_n)$ is a direct sum of copies of~$\cO_\dX$. Moreover, the coboundary maps marked~$\delta$ are isomorphisms, since the extension~\eqref{eq:RH_universal_extension} was chosen so that its coboundary map is the identity. Thus the central horizontal maps are both zero, implying that the third vertical map is injective as desired. This completes the inductive proof of Proposition~\ref{prop:RH_hodge_filtration}.
\end{proof}

The filtrations on~$\dE^\RH_n$ from Proposition~\ref{prop:RH_hodge_filtration} make~$\dE^\RH$ into a pro-object of~$\FMIC^\un(\dX,\cO_\dX)$. This filtration on~$\dE^\RH$ also admits an abstract characterisation via a universal property.

\begin{prop}\label{prop:RH_filtered_representability}
	The pair~$(\dE^\RH,e^\RH)$ pro-represents the functor
	\begin{align*}
		\FMIC^\un(\dX,\cO_\dX) &\to \Mod(\sB_\deR^+) \\
		\dE &\mapsto \rF^0\!\dE_x
	\end{align*}
	\begin{proof}
		Since~$(\dE^\RH,e^\RH)$ represents the fibre functor
		\[
		\omega^\RH_x\colon \MIC^\un(\dX,\cO_\dX) \to \Vec(\sB_\deR) \,,
		\]
		the proposition amounts to the following: for any~$\dE\in\FMIC^\un(\dX,\cO_\dX)$, a morphism $f\colon\dE^\RH\to\dE$ in $\pro{-}\MIC^\un(\dX,\cO_\dX)$ is filtered if and only if~$f(e^\RH)\in\rF^0\!\dE_x$. The only if direction is immediate since we have~$e^\RH\in\rF^0\!\dE^\RH_x$; we prove the converse implication by induction on~$\rank(\dE)$.
		
		In the base case $\dE=\cO_\dX$, since $\rH^0(\dX,\cO_\dX)=\sB_\deR$ we know from Proposition~\ref{prop:tannakian_explicit_universal} that any map $f\colon\dE^\RH\to\dE$ factors uniquely as
		\[
		\dE^\RH \to \dE^\RH_0 = \cO_\dX \xrightarrow{\alpha} \cO_\dX = \dE
		\]
		where the latter map is multiplication by some~$\alpha\in\sB_\deR$. So~$f$ is filtered if and only if~$\alpha\in\rF^0\!\sB_\deR$, which occurs if and only if $f(e^\RH)\in\rF^0\!\dE_x$ as desired.
		
		For the inductive step, we may write any non-zero~$\dE\in\FMIC^\un(\dX,\cO_\dX)$ as a strict extension
		\[
		0 \to \cO_\dX \to \dE \to \dE' \to 0 \,,
		\]
		and we may assume inductively that the claim holds for~$\dE'$. Fix a map~$f\colon\dE^\RH\to\dE$, and suppose that~$f(e^\RH)\in\rF^0\!\dE_x$. Writing~$f'$ for the composite
		\[
		\dE^\RH \xrightarrow{f} \dE \to \dE' \,,
		\]
		we have that~$f'(e^\RH)\in\rF^0\!\dE'_x$ and so~$f'$ is a filtered map by inductive assumption.
		
		Now~$f'$ factors through a filtered morphism~$f'_n\colon\dE^\RH_n\to\dE'$ for some~$n$. We saw in the proof of Proposition~\ref{prop:RH_hodge_filtration} that the map
		\[
		\rH^1_\FdeR(\dX,\dE^{\RH,\vee}_n) \to \rH^1_\FdeR(\dX,\dE^{\RH,\vee}_{n+1})
		\]
		is zero, so by Lemma~\ref{lem:filtered_RH_extensions} any strict extension of~$\dE^\RH_n$ by~$\cO_\dX$ splits when pulled back to~$\dE^\RH_{n+1}$. This implies that~$f'_n$ lifts to a filtered map~$f''_{n+1}\colon\dE^\RH_{n+1}\to\dE$, which we may further lift to a filtered map~$f''\colon\dE^\RH\to\dE$.
		
		Now the difference~$f-f''\colon\dE^\RH\to\cO_\dX$ satisfies~$f(e^\RH)-f''(e^\RH)\in\rF^0\!\sB_\deR$, so by the base case again we find that~$f-f''$ is a filtered map. Hence~$f=(f-f'')+f''$ is filtered, as desired.
	\end{proof}
\end{prop}

For any two points~$x,y\in Y^\an(L)$ defined over a finite extension~$L$ of~$K$ inside~$\bC_K$, we define the \emph{Hodge filtration} on $\pi_1^\RH(\dX;x,y)$ to be the filtration on its affine ring induced from the filtration on~${}_x\dE^\RH$ via the identification
\[
\cO(\pi_1^\RH(\dX;x,y)) \cong {}_x\dE^{\RH,\vee}_y
\]
of Proposition~\ref{prop:tannakian_groupoid_via_universal_objects}. This is an exhaustive, separated $\sB_\deR^+$-linear filtration. It follows from Proposition~\ref{prop:RH_filtered_representability} that the Hodge filtration on~${}_x\dE^\RH$ is compatible with the coalgebra structure, and hence the Hodge filtration on~$\pi_1^\RH(\dX;x,y)$ is compatible with the algebra structure on its affine ring. Moreover, the Hodge filtration is also compatible with identities, composition and path-reversal in the fundamental groupoid, which is a consequence of the following lemma, which also gives another characterisation of the Hodge filtration.

\begin{lem}\label{lem:RH_filtered_yoneda}
	The isomorphism~\eqref{eq:tannakian_yoneda}
	\[
	\Hom(\omega^\RH_x,\omega^\RH_y) \cong {}_x\dE^\RH_y
	\]
	is strictly compatible with filtrations, where the filtration on the right-hand side is the fibre of the Hodge filtration and the filtration on the left-hand is the one where a natural transformation~$\gamma\colon\omega^\RH_x\to\omega^\RH_y$ lies in~$\rF^i$ just when~$\gamma_\dE$ is filtered of degree~$i$ for all~$\dE\in\FMIC^\un(\dX,\cO_\dX)$.
	\begin{proof}
		One direction is easy: if~$\gamma\in\rF^i$ then~$\gamma(e^\RH_x)\in\rF^i\!{}_x\dE^\RH_y$. Conversely, if~$\gamma\in\Hom(\omega^\RH_x,\omega^\RH_y)$ corresponds to an element~$e'\in\rF^i\!{}_x\dE^\RH_y$, then for any~$\dE\in\MIC^\un(\dX,\cO_\dX)$ the map~$\gamma_\dE$ is given by
		\[
		\gamma_\dE(e) = f_e(e')
		\]
		where~$f_e\colon{}_x\dE^\RH\to\dE$ is the unique map such that~$f_e(e^\RH)=e$. If~$\dE$ is filtered and~$e\in\rF^j\!\dE_x$, then by Proposition~\ref{prop:RH_filtered_representability} applied to a shift of the filtration on~$\dE$ we see that~$f_e$ is filtered of degree~$j$, and hence~$\gamma_\dE(e)\in\rF^{i+j}\!\dE_y$. Hence~$\gamma$ is in~$\rF^i$ as desired.
	\end{proof}
\end{lem}
\section{Comparison between \'etale and de Rham fundamental groupoids}

Having defined the Riemann--Hilbert fundamental groupoid, we can now relate it to both the \'etale and de Rham fundamental groupoids, proving the comparison isomorphisms of Theorem~\ref{thm:comparison}.

\subsection{Comparison with de Rham fundamental groupoid}\label{ss:de_rham}

Let
\[
\MIC^\un(X,\cO_X) \hspace{0.4cm},\hspace{0.4cm} \MIC^\un(X^\an,\cO_{X^\an}) \hspace{0.4cm},\hspace{0.4cm} \MIC^\un(Y,\cO_Y) \hspace{0.4cm}\text{and}\hspace{0.4cm} \MIC^\un(Y^\an,\cO_{Y^\an})
\]
denote the categories of unipotent vector bundles with integrable (logarithmic) connection on~$X_\Zar$, $X^\an_\an$, $Y_\Zar$ and~$Y^\an_\an$, respectively, where~$(-)_\Zar$ and~$(-)_\an$ denote the Zariski and analytic sites.

These categories are all finitely generated unipotent Tannakian categories over~$K$, and the obvious pullback morphisms of sites induce~$\otimes$-functors
\begin{center}
\begin{tikzcd}
	\MIC^\un(X,\cO_X) \arrow[r,"(-)^\an"]\arrow[d,"(-)|_Y"] & \MIC^\un(X^\an,\cO_{X^\an}) \arrow[d,"(-)|_{Y^\an}"] \\
	\MIC^\un(Y,\cO_Y) \arrow[r,"(-)^\an"] & \MIC^\un(Y^\an,\cO_{Y^\an})
\end{tikzcd}
\end{center}
which are all equivalences. (See, variously, \cite[\S2.2]{deligne:local_behaviour} and~\cite[Theorem~9.1]{scholze:relative}.) Any $K$-rational point~$x\in Y(K)$ gives rise to compatible fibre functors on each of these four categories, which we denote by~$\omega_x^\deR$. We denote the corresponding Tannaka groupoids by
\[
\pi_1^\deR(X;-,-) \hspace{0.4cm},\hspace{0.4cm} \pi_1^\deR(X^\an;-,-) \hspace{0.4cm},\hspace{0.4cm} \pi_1^\deR(Y;-,-) \hspace{0.4cm}\text{and}\hspace{0.4cm} \pi_1^\deR(Y^\an;-,-) \,,
\]
respectively. For the sake of brevity we will confine our further discussion to the second of these only.

The de Rham fundamental groupoid~$\pi_1^\deR(X^\an;-,-)$ carries a Hodge filtration described in the case of affine algebraic curves in \cite[\S3]{hadian:motivic_pi_1}, which we now describe. For this, we let~$\FMIC^\un(X^\an,\cO_{X^\an})$ denote the category of \emph{unipotent filtered vector bundles with integrable logarithmic connection}, defined as follows.

\begin{defi}
	Let~$\dE$ be an~$\cO_{X^\an}$-module, $\nabla\colon\dE\to\dE\otimes_{\cO_{X^\an}}\Omega^1_{X^\an/K}$ an integrable logarithmic connection on~$\dE$, and~$\rF^\bullet$ a decreasing $\cO_{X^\an}$-linear filtration on~$\dE$ satisfying Griffiths transversality with respect to~$\nabla$. We say that~$(\dE,\nabla,\rF^\bullet)$ is \emph{unipotent} just when there is a finite $\nabla$-stable $\cO_{X^\an}$-linear filtration on~$\dE$, each of whose graded pieces is filtered-isomorphic to
	\[
	V\otimes_K\cO_{X^\an}
	\]
	for some finitely filtered finite-dimensional $K$-vector space~$V$. We denote the category of such triples~$(\dE,\nabla,\rF^\bullet)$ by~$\FMIC^\un(X^\an,\cO_{X^\an})$.
\end{defi}

\begin{rmk}
	The condition on the triple~$(\dE,\nabla,\rF^\bullet)$ in the above definition is stronger than merely requiring~$(\dE,\nabla)$ to be unipotent. In fact (though we will not use this), if~$(\dE,\nabla)$ is unipotent and~$\rF^\bullet$ is a Griffiths-transverse filtration on~$\dE$, then~$(\dE,\nabla,\rF^\bullet)\in\FMIC^\un(X^\an,\cO_{X^\an})$ if and only if~$\rF^\bullet$ is a filtration by \emph{unipotent} $\cO_{X^\an}$-submodules.
	
\end{rmk}

Suppose now that~$x\in Y^\an(K)$ is a $K$-rational point. Let~$(\dE^\deR,e^\deR) = ({}_x\dE^\deR,e_x^\deR)$ be the pair pro-representing~$\omega_x^\deR$ as in \S\ref{ss:tannakian_universal_objects}, and write~$\dE^\deR$ as the inverse limit of~$(\dE^\deR_n)_{n\geq0}$ as in Proposition~\ref{prop:tannakian_explicit_universal}. The Hodge filtration on~$\dE^\deR$ is then characterised recursively as follows.

\begin{prop}[{cf.\ \cite[Lemma~3.6]{hadian:motivic_pi_1}}]\label{prop:dR_hodge_filtration}
	There is a unique way to put a decreasing filtration~$\rF^\bullet$ on each $\dE^\deR_n$ so that the following properties are satisfied:
	\begin{enumerate}\setcounter{enumi}{-1}
		\item\label{condn:dR_hodge_basic} each $\dE^\deR_n$ is an object of $\FMIC^\un(X^\an,\cO_{X^\an})$;
		\item\label{condn:dR_hodge_normalisation} the filtration on~$\dE^\deR_0=\cO_{X^\an}$ is the trivial filtration;
		\item\label{condn:dR_hodge_exactness} for all~$n\geq1$ the sequence
		\begin{equation}\label{eq:dR_universal_extension}
			0 \to \cO_{X^\an}\otimes_K\rH^1_\deR(X^\an,\dE^{\deR,\vee}_{n-1})^\vee \to \dE^\deR_n \to \dE^\deR_{n-1} \to 0
		\end{equation}
		is strict exact for the filtrations; and
		\item\label{condn:dR_hodge_rigidity} $e^\deR_n\in\rF^0\!\dE^\deR_{n,x}$ for all~$n$.
	\end{enumerate}
	The sequence~\eqref{eq:dR_universal_extension} in the second point is the sequence~\eqref{eq:tannakian_universal_extension}, using the usual identification of Ext-groups in~$\MIC^\un(X^\an,\cO_{X^\an})$ with de Rham cohomology. The implicit filtration on~$\rH^1_\deR(X^\an,\dE^{\deR,\vee}_{n-1})$ is the natural one induced from the filtration on~$\dE^\deR_{n-1}$.
\end{prop}

\begin{rmk}
	The corresponding statement of \cite[Lemma~3.6]{hadian:motivic_pi_1} omits condition~\eqref{condn:dR_hodge_rigidity}: without this condition the filtration is not unique, only unique up to the action of~$\Aut(\dE^\deR_n)$.
\end{rmk}

\begin{proof}
	A proof similar to the proof of Proposition~\ref{prop:RH_hodge_filtration} works; we give just the outline. One first checks, as in Lemma~\ref{lem:filtered_RH_extensions}, that the groups of strict extensions in~$\FMIC^\un(X^\an,\cO_{X^\an})$ are computed by
	\[
	\Ext^1_{\FMIC}(\dE_1,\dE_2) \cong \rH^1_\FdeR(X^\an,\dE_1^\vee\otimes_{\cO_{X^\an}}\dE_2) \colonequals \bH^1(X^\an_\an,\rF^0\!\DeR(\dE)) \,.
	\]
	Using this, one shows inductively that there exist unique filtrations on the~$\dE^\deR_n$ satisfying~\eqref{condn:dR_hodge_basic}--\eqref{condn:dR_hodge_rigidity}, and that these filtrations additionally satisfy the auxiliary condition
	\textit{
	\begin{enumerate}\setcounter{enumi}{3}
		\item\label{condn:dR_hodge_good} the natural map $\rH^1_\FdeR(X^\an,\dE^\deR_n\langle i\rangle)\to\rH^1_\deR(X^\an,\dE^\deR_n)$ is injective for all~$i$, where~$(-)\langle i\rangle$ denotes a shift of the filtration by~$i$.
	\end{enumerate}
	}
	In the base case~$\dE^\deR_0=\cO_{X^\an}$, condition~\eqref{condn:dR_hodge_good} follows from the degeneration of the Hodge--de Rham spectral sequence at the first page. From here on, the induction runs similarly to in Proposition~\ref{prop:RH_hodge_filtration}.
\end{proof}

The filtration on~${}_x\dE^\deR=\dE^\deR=\varprojlim_n\dE^\deR_n$ provided by Proposition~\ref{prop:dR_hodge_filtration} is called the \emph{Hodge filtration}, and induces a Hodge filtration on the affine ring of $\pi_1^\deR(X^\an;x,y)$ for all~$x,y\in Y^\an(K)$ via Proposition~\ref{prop:tannakian_groupoid_via_universal_objects}. Arguing as in \S\ref{ss:RH_hodge_filtration}, this Hodge filtration is compatible with the algebra structure on the affine ring, and with identities, composition and inversion in the fundamental groupoid.
\smallskip

To compare the de Rham and Riemann--Hilbert fundamental groupoids, we will consider the functor $\MIC^\un(X^\an,\cO_{X^\an})\to\MIC^\un(\dX,\cO_\dX)$ given by
\[
\dE \mapsto \sB_\deR\hatotimes_K\dE = \cO_\dX\otimes_{\cO_{X^\an}}\dE \,.
\]
This is a $K$-linear exact $\otimes$-functor, and the induced map on first Ext-groups is given by the induced map on de Rham cohomology. There is an induced $G_K$-action on~$\cO_\dX\otimes_{\cO_{X^\an}}\dE$, and, for~$\dE\in\FMIC^\un(X^\an,\cO_{X^\an})$, also an induced filtration in the sense of Definition~\ref{def:RH_filtrations}. Regarding the Ext-groups, we have the following.

\begin{lem}\label{lem:RH-dR_extension_comparison}
	For every~$\dE\in\MIC^\un(X^\an,\cO_{X^\an})$, the induced map
	\[
	\sB_\deR\otimes_K\rH^1_\deR(X^\an,\dE) \to \rH^1_\deR(\dX,\sB_\deR\hatotimes_K\dE)
	\]
	is a $G_K$-equivariant isomorphism for all~$i$. If~$\dE\in\FMIC^\un(X^\an,\cO_{X^\an})$, it is moreover strictly compatible with filtrations.
	\begin{proof}
		The map in question is clearly $G_K$-equivariant. For the remaining parts we make the stronger claim that the map
		\begin{equation}\label{eq:RH-dR_derived_comparison}\tag{$\ast$}
		\sB_\deR\otimes_K\rR\!\Gamma(X^\an_\an,\DeR(\dE)) \to \rR\!\Gamma(\dX,\DeR(\sB_\deR\hatotimes_K\dE))
		\end{equation}
		is a filtered quasi-isomorphism for all~$\dE\in\FMIC^\un(X^\an,\cO_{X^\an})$. Indeed, it is an isomorphism for~$\dE=\cO_{X^\an}$ by Lemma~\ref{lem:RH-dR_extension_comparison_baby}, and hence also for~$\dE$ any filtration-shift of~$\cO_{X^\an}$. Since the property of~\eqref{eq:RH-dR_derived_comparison} being a filtered quasi-isomorphism is preserved under strict extensions, it is a filtered quasi-isomorphism for every~$\dE\in\FMIC^\un(X^\an,\cO_{X^\an})$ and we are done.
	\end{proof}
\end{lem}

By the general machinery of \S\ref{s:tannakian}, Lemma~\ref{lem:RH-dR_extension_comparison} translates into a comparison result between the unipotent de Rham and Riemann--Hilbert fundamental groupoids. For any~$x\in Y^\an(K)$, the square
\begin{center}
\begin{tikzcd}[column sep = large]
	\MIC^\un(X^\an,\cO_{X^\an}) \arrow[r,"\sB_\deR\hatotimes_K(-)"]\arrow[d,"\omega_x^\deR"] & \MIC^\un(\dX,\cO_\dX) \arrow[d,"\omega_x^\RH"] \\
	\Vec(K) \arrow[r,"\sB_\deR\otimes_K(-)"] & \Vec(\sB_\deR)
\end{tikzcd}
\end{center}
commutes (up to a canonical $\otimes$-natural isomorphism). So by the discussion in \S\ref{ss:tannakian_comparison} there is a natural comparison map
\begin{equation}\label{eq:RH-dR_functoriality_universal}
	{}_x\dE^\RH \to \sB_\deR\hatotimes_K{}_x\dE^\deR
\end{equation}
in $\pro{-}\MIC^\un(\dX,\cO_\dX)$, and for any second $K$-rational point~$y\in Y^\an(K)$ there is also a comparison map
\begin{equation}\label{eq:RH-dR_functoriality}
	\pi_1^\RH(\dX;x,y) \to \sB_\deR\otimes_K\pi_1^\deR(X^\an;x,y)
\end{equation}
in the category of~$\sB_\deR$-schemes. These maps are compatible with all groupoid operations. We have thus proved the following.

\begin{prop}\label{prop:dR-RH_comparison}
	The maps~\eqref{eq:RH-dR_functoriality_universal} and~\eqref{eq:RH-dR_functoriality} are isomorphisms for all points~$x,y\in Y^\an(K)$, and are $G_K$-equivariant and strictly compatible with Hodge filtrations.
	\begin{proof}
		Using Lemma~\ref{lem:RH-dR_extension_comparison}, Proposition~\ref{prop:tannakian_iso} shows that the induced maps are isomorphisms. For $G_K$-equivariance, we note that~\eqref{eq:RH-dR_functoriality_universal} takes~$e^\RH_x$ to~$1\otimes e^\deR_x\in(\sB_\deR\otimes_K{}_x\dE^\deR_x)^{G_K}$. This implies by Lemma~\ref{lem:RH_equivariant_representability} that~\eqref{eq:RH-dR_functoriality_universal} is $G_K$-equivariant, and this automatically implies the same for~\eqref{eq:RH-dR_functoriality}.
		
		The argument for Hodge-filteredness is slightly more subtle. It suffices to show that the isomorphisms
		\[
		\phi_n\colon {}_x\dE^\RH_n\to\sB_\deR\hatotimes_K{}_x\dE^\deR
		\]
		of~\eqref{eq:tannakian_functoriality_family_levelwise} are filtered isomorphisms for all~$n$. Proposition~\ref{prop:RH_filtered_representability} shows that~$\phi_n$ is a filtered map, but it is not immediately obvious that it is strict. We prove this by induction, the base case~$n=0$ being trivial. Supposing inductively that~$\phi_n$ is strict for the Hodge filtration, we recall from the proof of Proposition~\ref{prop:tannakian_iso} that we have a commuting diagram
		\begin{center}
		\begin{tikzcd}
			0 \arrow[r] & \cO_\dX\otimes_{\sB_\deR}\rH^1_\deR(\dX,{}_x\dE^{\RH,\vee}_n)^\vee \arrow[r]\arrow[d,"\wr"'] & {}_x\dE^\RH_{n+1} \arrow[r]\arrow[d,"\phi_{n+1}","\wr"'] & {}_x\dE^\RH_n \arrow[r]\arrow[d,,"\phi_n","\wr"'] & 0 \\
			0 \arrow[r] & \cO_\dX\otimes_K\rH^1_\deR(X^\an,{}_x\dE^{\deR,\vee}_n)^\vee \arrow[r] & \sB_\deR\hatotimes_K{}_x\dE^\deR_{n+1} \arrow[r] & \sB_\deR\hatotimes_K{}_x\dE^\deR_n \arrow[r] & 0
		\end{tikzcd}
		\end{center}
		with strict exact rows, where the unmarked vertical map is dual to the composite map
		\[
		\sB_\deR\otimes_K\rH^1_\deR(X^\an,{}_x\dE^{\deR,\vee}_n) \xrightarrow{\sim} \rH^1_\deR(\dX,\sB_\deR\hatotimes{}_x\dE^{\deR,\vee}_n) \xrightarrow[\sim]{\phi_n^*}  \rH^1_\deR(\dX,{}_x\dE^{\RH,\vee}_n) \,.
		\]
		By Lemma~\ref{lem:RH-dR_extension_comparison} this map is strict for the filtrations. So by the five-lemma~$\phi_{n+1}$ is also strict, completing the induction.
	\end{proof}
\end{prop}

\subsection{Comparison with \'etale fundamental groupoid}

The comparison between \'etale and Riemann--Hilbert fundamental groups is a little more involved than the de Rham case. Let
\[
\Loc^\un(X_{\bC_K},\bQ_p) \hspace{0.4cm},\hspace{0.4cm} \Loc^\un(X^\an_{\bC_K},\bQ_p) \hspace{0.4cm},\hspace{0.4cm} \Loc^\un(Y_{\bC_K},\bQ_p) \hspace{0.4cm}\text{and}\hspace{0.4cm} \Loc^\un(Y^\an_{\bC_K},\bQ_p)
\]
denote the categories of unipotent $\bQ_p$-local systems on~$X_{\bC_K,\ket}$, $X^\an_{\bC_K,\ket}$, $Y_{\bC_K,\et}$ and~$Y^\an_{\bC_K,\et}$, respectively. Here~$(-)_\ket$ denotes the Kummer \'etale site, see \cite[\S2]{illusie} and \cite[\S4.1]{dllz:foundations}, respectively.

These categories are all finitely generated unipotent Tannakian categories over~$\bQ_p$, and the obvious pullback morphisms of sites induce $\otimes$-functors
\begin{center}
\begin{tikzcd}
	\Loc^\un(X_{\bC_K},\bQ_p) \arrow[r,"(-)^\an"]\arrow[d,"(-)|_Y"] & \Loc^\un(X^\an_{\bC_K},\bQ_p) \arrow[d,"(-)|_{Y^\an}"] \\
	\Loc^\un(Y_{\bC_K},\bQ_p) \arrow[r,"(-)^\an"] & \Loc^\un(Y^\an_{\bC_K},\bQ_p) \,,
\end{tikzcd}
\end{center}
which are all equivalences. (See, variously, \cite[Example~4.7(c)]{illusie}, \cite[Theorem~9.3]{scholze:relative} and \cite[Corollary~6.3.4]{dllz:foundations}.) Any point~$x\in Y(\bC_K)$ gives rise to compatible fibre functors on each of these four categories, which we denote by~$\omega_{\bar x}^\ket$ or~$\omega_{\bar x}^\et$. We denote the corresponding Tannaka groupoids by
\[
\pi_1^\ket(X_{\bC_K};-,-) \hspace{0.4cm},\hspace{0.4cm} \pi_1^\ket(X^\an_{\bC_K};-,-) \hspace{0.4cm},\hspace{0.4cm} \pi_1^\et(Y_{\bC_K};-,-) \hspace{0.4cm}\text{and}\hspace{0.4cm} \pi_1^\et(Y^\an_{\bC_K};-,-) \,,
\]
respectively. For the sake of brevity we will confine our further discussion to the first of these only, but a similar discussion holds in all cases.

There is a canonical left\footnote{Actually, the action is more natural from the right.} action of~$G_K$ on the site~$X_{\bC_K,\ket}$, and hence a natural right action on $\Loc^\un(X_{\bC_K},\bQ_p)$. If~$\bar x=x\in Y(K)$ is a $K$-rational point, then its fibre functor $\omega_{\bar x}^\ket\colon\Loc^\un(X_{\bC_K},\bQ_p)\to\Vec(\bQ_p)$ is $G_K$-invariant. Hence for any second $K$-rational point $y\in Y(K)$, the functor of $\otimes$-natural isomorphisms $\omega_{\bar x}^\ket\xrightarrow\sim\omega_{\bar y}^\ket$ is $G_K$-invariant, and hence there is an induced left $G_K$-action on~$\pi_1^\ket(X_{\bC_K};\bar x,\bar y)$. For varying~$\bar x$ and~$\bar y$, these actions clearly respect the groupoid structure maps.

Similarly, there is an action of~$G_K$ on the universal object ${}_{\bar x}\bar\bE^\ket$ lying over the action of~$G_K$ on~$X_{\bC_K}$, i.e.\ there are isomorphisms
\[
\phi_\sigma\colon {}_{\bar x}\bar\bE^\ket \to \sigma^*({}_{\bar x}\bar\bE^\ket)
\]
namely the unique map sending $e_{\bar x}^\ket\in{}_{\bar x}\bar\bE^\ket_{\bar x}$ to $e_{\bar x}^\ket\in (\sigma^*({}_{\bar x}\bar\bE^\ket))_{\bar x} = {}_{\bar x}\bar\bE^\ket_{\bar x}$. This action can also be described by universal properties similar to Lemmas~\ref{lem:RH_equivariant_representability} and~\ref{lem:RH_yoneda_equivariant}, and hence the $G_K$-actions we have described on $\pi_1^\ket(X_{\bC_K};-,-)$ and on~${}_{\bar x}\bar\bE^\ket$ agree via the identification in Proposition~\ref{prop:tannakian_groupoid_via_universal_objects}.
\smallskip

The isomorphisms~$\phi_\sigma$ above satisfy the cocycle condition
\begin{equation}\label{eq:cocycle_condition}
\phi_{\sigma\tau} = \tau^*(\phi_\sigma)\circ\phi_\tau \,,
\end{equation}
and so can be viewed as descent data for the ``covering'' $X_{\bC_K}\to X$. An important technical point for us is that this descent data is effective.

\begin{lem}\label{lem:descent_of_universal}
	There exists a pro-$\bQ_p$-local system~${}_x\bE^\ket$ on~$X_\ket$ whose pullback to~$X_{\bC_K,\ket}$ is isomorphic to~${}_{\bar x}\bar\bE^\ket$ and such that the isomorphisms~$\phi_\sigma$ are the natural ones.
\end{lem}

We will deduce this as a special case of a more general criterion. For this, we will say that a $\bQ_p$-local system~$\bE$ on~$X_\ket$ is \emph{geometrically unipotent} just when it has a finite filtration, all of whose graded pieces are pullbacks of $\bQ_p$-local systems on~$\Spec(K)_\et$.

\begin{lem}
	Let~$\bar\bE\in\Loc^\un(X_{\bC_K},\bQ_p)$, and suppose we are given isomorphisms $\phi_\sigma\colon\bar\bE\xrightarrow\sim\sigma^*(\bar\bE)$ for all~$\sigma\in G_K$ satisfying the cocycle condition~\eqref{eq:cocycle_condition}. Then~$(\bar\bE,(\phi_\sigma)_{\sigma\in G_K})$ is the pullback of a geometrically unipotent $\bQ_p$-local system~$\bE$ on~$X_\ket$ if and only if the induced $G_K$-action on~$\bar\bE_{\bar x}$ is continuous. Moreover,~$\bE$ is unique up to unique isomorphism.
	\begin{proof}
		Note that every unipotent $\bQ_p$-local system on~$X_{\bC_K,\ket}$ is an isogeny $\bZ_p$-local system, meaning that there is a (unipotent) $\bZ_p$-local system~$\bar\bE_0$ such that~$\bar\bE=\bar\bE_0[p^{-1}]$. Hence $\Loc^\un(X_{\bC_K},\bQ_p)$ is equivalent via the fibre functor $\omega_{\bar x}^\ket$ to the category~$\Rep_\cts^\un(\hat\pi_1^\ket(X_{\bC_K},\bar x))$ of unipotent continuous representations of the profinite Kummer \'etale fundamental group $\hat\pi_1^\ket(X_{\bC_K},\bar x)$~\cite[\S4.6]{illusie}. Similarly, the category~$\Loc^\gun(X,\bQ_p)$ of geometrically unipotent $\bQ_p$-local systems on~$X$ is equivalent to the category~$\Rep_\cts^\gun(\hat\pi_1^\ket(X,\bar x))$ of continuous representations of~$\hat\pi_1^\ket(X,\bar x)=\hat\pi_1^\ket(X_{\bC_K},\bar x)\rtimes G_K$ which possess a filtration such that the action on each graded piece factors through~$G_K$. Under these correspondences, the pullback functor~$\Loc^\gun(X,\bQ_p)\to\Loc^\un(X_{\bC_K},\bQ_p)$ is the forgetful functor
		\[
		\Rep_\cts^\gun(\hat\pi_1^\ket(X,\bar x)) \to \Rep_\cts^\un(\hat\pi_1^\ket(X_{\bC_K},\bar x)) \,.
		\]
		
		Giving isomorphisms~$(\phi_\sigma)_{\sigma\in G_K}$ on some~$\bar\bE\in\Loc^\un(X_{\bC_K},\bQ_p)$ satisfying~\eqref{eq:cocycle_condition} is equivalent to giving an action of~$G_K$ on an object~$V\in\Rep_\cts^\un(\hat\pi_1^\ket(X_{\bC_K},\bar x))$ satisfying
		\[
		\sigma(\gamma\cdot v) = \sigma(\gamma)\cdot\sigma(v)
		\]
		for all~$\sigma\in G_K$, $\gamma\in\hat\pi_1^\ket(X_{\bC_K},\bar x)$ and~$v\in V$. In other words, it corresponds to giving an action of~$\hat\pi_1^\ket(X,\bar x)$ on~$V$ extending the given action of~$\hat\pi_1^\ket(X_{\bC_K},\bar x)$. It is then clear that~$V$ comes from an object $V_0\in\Rep_\cts^\gun(\hat\pi_1^\ket(X,\bar x))$ if and only if the $G_K$-action on~$V$ is continuous, in which~$V_0$ is unique up to isomorphism. This is what we wanted to show.
	\end{proof}
\end{lem}

\begin{proof}[Proof of Lemma~\ref{lem:descent_of_universal}]
	$\pi_1^\ket(X_{\bC_K},\bar x)$ is the topological Mal\u cev completion of~$\hat\pi_1^\ket(X_{\bC_K},\bar x)$ and so we have
	\[
	{}_{\bar x}\bar\bE^\ket_{\bar x} = \cO(\pi_1^\ket(X_{\bC_K},\bar x))^\vee = \bQ_p\llbrack\hat\pi_1^\ket(X_{\bC_K},\bar x)\rrbrack^\wedge \,,
	\]
	where~$(-)^\wedge$ denotes completion with respect to the augmentation ideal. These identifications are compatible with $G_K$-actions, and the action on the latter is certainly pro-continuous.
\end{proof}

In preparation to compare the pro-unipotent Kummer \'etale fundamental groupoid with the pro-unipotent Riemann--Hilbert fundamental groupoid, let us give an equivalent description of the Tannakian category $\Loc^\un(X^\an_{\bC_K},\bQ_p)$ in terms of the pro-Kummer \'etale site \cite[\S5.1]{dllz:foundations}. We let~$X^\an_\proket$ and~$X^\an_{\bC_K,\proket}$ denote the pro-Kummer \'etale sites of~$X^\an$ and~$X^\an_{\bC_K}$, respectively, which fit into a 2-commuting diagram of sites
\begin{center}
	\begin{tikzcd}
		X^\an_{\bC_K,\ket} \arrow[d,"\beta_\ket"] & X^\an_{\bC_K,\proket} \arrow[l,"\bar\nu"']\arrow[r,"\bar\mu"]\arrow[d,"\beta_\proket"] & \dX \arrow[d,equals] \\
		X^\an_\ket & X^\an_\proket \arrow[l,"\nu"']\arrow[r,"\mu"] & X^\an_\an \,.
	\end{tikzcd}
\end{center}

As in~\cite[Definition~6.3.2]{dllz:foundations} we write~$\hat\bZ_p\colonequals\varprojlim_n(\bZ/p^n)$ for the inverse limit of the constant sheaves~$\bZ/p^n$ on~$X^\an_\proket$, and set~$\hat\bQ_p\colonequals\hat\bZ_p[p^{-1}]$. A \emph{$\hat\bQ_p$-local system} on~$X^\an_\proket$ or~$X^\an_{\bC_K,\proket}$ means a locally finite free sheaf of~$\hat\bQ_p$-modules. We write~$\Loc^\un(X^\an_{\bC_K},\hat\bQ_p)$ for the category of unipotent $\hat\bQ_p$-local systems on~$X^\an_{\bC_K,\proket}$. Any~$\bC_K$-point $\bar x\in Y^\an(\bC_K)$ determines a morphism of sites
\[
i_{\bar x}\colon\Sp(\bC_K)_\proet \to X^\an_{\bC_K,\proket}
\]
where~$\Sp(\bC_K)_\proet$ is the pro-\'etale site of \cite[Definition~3.9]{scholze:relative}\cite[Erratum~1]{scholze:erratum} (it is also the pro-Kummer \'etale site of~$\Sp(\bC_K)$ for the trivial log structure). We define the functor
\begin{align*}
	\omega^\proket_{\bar x}\colon\Loc^\un(X_{\bC_K},\hat\bQ_p) &\to \Mod(\bQ_p) \\
	\hat{\bar\bE} &\mapsto \rH^0(\Sp(\bC_K),i_{\bar x}^{-1}\hat{\bar\bE}) \,.
\end{align*}

\begin{prop}
	$\Loc^\un(X^\an_{\bC_K},\hat\bQ_p)$ is a $\bQ_p$-linear finitely generated unipotent neutral Tannakian category, and~$\omega^\proket_{\bar x}$ is a fibre functor for all~$\bar x\in Y^\an(\bC_K)$.
	\begin{proof}
		We verify the conditions of Lemma~\ref{lem:tannakian_criteria}. Only two of the conditions require discussion: finite-dimensionality of~$\Ext^1(\hat\bQ_p,\hat\bQ_p)$ and exactness of~$\omega^\proket_{\bar x}$. The former will follow from Lemma~\ref{lem:proetale_local_systems} shortly, for the latter it suffices to show that every pro-\'etale covering of~$\Sp(\bC_K)$ splits. This follows from the definition of the pro-\'etale site: any pro-\'etale map $U\to\Sp(\bC_K)$ factors as $U\to U_0\to\Sp(\bC_K)$ for~$U_0\to\Sp(\bC_K)$ \'etale (so a disjoint union of copies of~$\Sp(\bC_K)$) and $U\to U_0$ a cofiltered inverse limit of finite \'etale maps. As soon as $U\neq\emptyset$, we have that~$U_0\to\Sp(\bC_K)$ splits, and so does $U\to U_0$ since a cofiltered inverse limit of finite non-empty sets is non-empty.
	\end{proof}
\end{prop}

Given a unipotent $\bQ_p$-local system~$\bar\bE$ on~$X^\an_{\bC_K}$, we may choose a $\bZ_p$-lattice~$\bar\bE_0=\varprojlim_n\bar\bE_{0,n}$ in~$\bar\bE$, and define a unipotent $\hat\bQ_p$-local system $\hat{\bar\bE}$ by
\[
\hat{\bar\bE} \colonequals \left(\varprojlim_n\bar\nu^{-1}\bar\bE_{0,n}\right)[p^{-1}] \,.
\]
It is easy to see that~$\hat{\bar\bE}$ is independent of the choice of~$\bZ_p$-lattice~$\bE_0$, and the assignment~$\bar\bE\mapsto\hat{\bar\bE}$ defines a $\otimes$-functor $\Loc^\un(X^\an_{\bC_K},\bQ_p)\to\Loc^\un(X^\an_{\bC_K},\hat\bQ_p)$.

\begin{lem}\label{lem:proetale_local_systems}
	The above functor $\Loc^\un(X^\an_{\bC_K},\bQ_p)\to\Loc^\un(X^\an_{\bC_K},\hat\bQ_p)$ is a $\otimes$-equivalence.
	\begin{proof}
		It suffices to check that~$\bar\bE\mapsto\hat{\bar\bE}$ induces an isomorphism
		\[
		\rH^i(X^\an_{\bC_K,\ket},\bar\bE) \to \rH^i(X^\an_{\bC_K,\proket},\hat{\bar\bE})
		\]
		for all~$i$ and all~$\bar\bE\in\Loc^\un(X^\an_{\bC_K},\bQ_p)$. In the case $\bar\bE=\bQ_p$, this follows from \cite[Corollary~6.3.4]{dllz:foundations}; in general it follows by induction.
	\end{proof}
\end{lem}


Now the pro-Kummer \'etale site~$X^\an_\proket$ carries a natural period sheaf\footnote{This sheaf is denoted~$\cO\bB_{\deR,\log}$ in~\cite{dllz:riemann-hilbert}, and in the case of trivial log structure is the completion of Scholze's period sheaf \cite[Definition~6.8(iv)]{scholze:relative}\cite[Erratum~3]{scholze:erratum} with respect to its filtration, see \cite[Remark~2.2.11]{dllz:riemann-hilbert}.}~$\cO\bB_\deR^\wedge$ \cite[Definition~2.2.10]{dllz:riemann-hilbert}. This is a sheaf of filtered~$\mu^{-1}\cO_{X^\an}$-algebras, and comes with an integrable logarithmic connection
\[
\nabla\colon \cO\bB_\deR^\wedge \to \cO\bB_\deR^\wedge \otimes_{\mu^{-1}\cO_{X^\an}} \mu^{-1}\Omega^1_{X^\an/K}
\]
which is compatible with the $\mu^{-1}\cO_{X^\an}$-algebra structure and satisfies Griffiths transversality with respect to the filtration. We have $\bar\mu_*(\cO\bB_\deR^\wedge) = \sB_\deR\hatotimes_K\cO_{X^\an} = \cO_\dX$ \cite[Lemma~3.3.2]{dllz:riemann-hilbert}, and this identification is compatible with connections, $G_K$-actions and filtrations.

The period sheaf~$\cO\bB_\deR^\wedge$ is also a sheaf of~$\hat\bQ_p$-algebras, where~$\hat\bQ_p$ is contained in the kernel of the connection~$\nabla$. Diao, Lan, Liu and Zhu define a unipotent Riemann--Hilbert functor as follows.

\begin{defi}[{\cite[\S3.2]{dllz:riemann-hilbert}\footnote{Strictly speaking, this is not quite the same Riemann--Hilbert functor as defined in~\cite[\S3.2]{dllz:riemann-hilbert}, since the functor described there is defined on $\hat\bQ_p$-local systems on~$X^\an_\proket$, whereas our~$\dRH$ will be defined on unipotent $\hat\bQ_p$-local systems on~$X^\an_{\bC_K,\proket}$.}}]
	The \emph{unipotent Riemann--Hilbert functor}
	\[
	\dRH\colon \Loc^\un(X^\an_{\bC_K},\bQ_p) \to \MIC^\un(\dX,\cO_\dX)
	\]
	is defined by
	\[
	\dRH(\bar\bE) \colonequals \bar\mu_*(\cO\bB_\deR^\wedge\otimes_{\hat\bQ_p}\hat{\bar\bE}) \,.
	\]
	The connection on~$\dRH(\bar\bE)$ is the natural one induced from the connection on~$\cO\bB_\deR^\wedge$. The Hodge filtration on~$\cO\bB_\deR^\wedge$ induces a filtration on~$\dRH(\bar\bE)$ for all~$\bar\bE$. If~$\bar\bE$ is the pullback of a geometrically unipotent $\bQ_p$-local system~$\bE$ on~$X^\an_\ket$, then there is an induced $G_K$-action on~$\dRH(\bar\bE)$.
\end{defi}

Note that it is not immediately obvious that the image of~$\dRH$ actually lands in~$\MIC^\un(\dX,\cO_\dX)$, since we do not yet know that~$\dRH$ is exact. This is fixed by the following lemma.

\begin{lem}\label{lem:etale-RH_extension_comparison}
	The functor~$\dRH$ is a $\bQ_p$-linear exact $\otimes$-functor, and for every~$\bar\bE\in\Loc^\un(X^\an_{\bC_K},\bQ_p)$ the induced map
	\begin{equation}\label{eq:etale-RH-extension_comparison}\tag{$\ast$}
	\sB_\deR\otimes_{\bQ_p}\rH^i_\ket(X^\an_{\bC_K},\bar\bE) \to \rH^i_\deR(\dX,\dRH(\bar\bE))
	\end{equation}
	is a filtered isomorphism for all~$i$. If~$\bar\bE$ is the pullback of a geometrically unipotent $\bQ_p$-local system on~$X^\an_\ket$, then the isomorphism is moreover $G_K$-equivariant.
	\begin{proof}
		\cite[Theorem~3.2.3(3)]{dllz:riemann-hilbert} shows that~\eqref{eq:etale-RH-extension_comparison} is a $G_K$-equivariant filtered isomorphism whenever~$\bar\bE$ is a pullback of a geometrically unipotent $\bQ_p$-local system on~$X^\an_\ket$. In particular it is an isomorphism for~$\bar\bE=\bQ_p$ the unit object. Moreover, \cite[Lemma~3.3.2]{dllz:riemann-hilbert} shows that the natural map~$\dRH(\bQ_p)\to\cO_\dX$ is an isomorphism, and that
		\[
		\rR^i\!\bar\mu_*(\cO\bB_\deR^\wedge) = 0
		\]
		for~$i>0$. It then follows by induction that
		\[
		\rR^i\!\bar\mu_*(\cO\bB_\deR^\wedge\otimes_{\hat\bQ_p}\hat{\bar\bE}) = 0
		\]
		for all~$\bar\bE\in\Loc^\un(X^\an_{\bC_K},\bQ_p)$ and all~$i>0$. So~$\dRH$ is exact and its image is contained in~$\MIC^\un(\dX,\cO_\dX)$. Since~$\dRH$ is certainly a lax $\otimes$-functor, it is in fact a $\otimes$-functor using Lemma~\ref{lem:tannakian_criteria}\eqref{condn:tannakian_lax_fibre_functor}.
		
		A similar inductive argument establishes that~\eqref{eq:etale-RH-extension_comparison} is an isomorphism for all~$\bar\bE\in\Loc^\un(X^\an_{\bC_K},\bQ_p)$. Showing that it is filtered requires a modicum of care, since the property that~\eqref{eq:etale-RH-extension_comparison} is a \emph{filtered} isomorphism is not necessarily preserved under extensions. We prove instead the \emph{a priori} stronger claim that the natural maps
		\begin{equation}\label{eq:etale-RH-extension_comparison_filtered}\tag{$\ast\ast$}
		\sB_\deR^+\otimes_{\bQ_p}\rH^i_\ket(X^\an_{\bC_K},\bar\bE) \to \rH^i_\FdeR(\dX,\dRH(\bar\bE))
		\end{equation}
		are isomorphisms for all~$\bar\bE\in\Loc^\un(X^\an_{\bC_K},\bQ_p)$ and all~$i$. This does hold for~$\bar\bE=\bQ_p$ using the fact that~\eqref{eq:etale-RH-extension_comparison} is a filtered isomorphism and the maps $\rH^i_\FdeR(\dX,\cO_\dX)\to\rH^i_\deR(\dX,\cO_\dX)$ are injective by Lemma~\ref{lem:OB_deR_good}, and this stronger property is preserved under extensions.
	\end{proof}
\end{lem}

Now for any $K$-rational point~$x\in Y^\an(K)$ we have a commuting square of ringed sites
\begin{center}
\begin{tikzcd}
	(X^\an_{\bC_K,\proket},\cO\bB_\deR^\wedge) \arrow[r,"\bar\mu"] & (\dX,\cO_\dX) \\
	(\Sp(\bC_K)_\profet,\cO\bB_{\deR,\Sp(\bC_K)}^\wedge) \arrow[r]\arrow[u,"i_{\bar x}"'] & (\ast,\sB_\deR) \arrow[u,"i_x"'] \,.
\end{tikzcd}
\end{center}
The base-change map from this square gives a map
\[
\dRH(\bar\bE)_x \to \sB_\deR\otimes_{\bQ_p}\hat{\bar\bE}_{\bar x} = \sB_\deR\otimes_{\bQ_p}\bar\bE_{\bar x} \,,
\]
$\otimes$-natural in~$\bar\bE\in\Loc^\un(X^\an_{\bC_K},\bQ_p)$. Since this base change map is an isomorphism for~$\bar\bE=\bQ_p$, it is an isomorphism in general, so the fibre functors~$\omega^\ket_{\bar x}$ and~$\omega^\RH_x$ are associated in the sense of~\S\ref{ss:tannakian_comparison}. This implies that we have a natural comparison map
\begin{equation}\label{eq:RH-et_functoriality_universal}
	{}_x\dE^\RH \to \dRH({}_{\bar x}\bar\bE^\ket)
\end{equation}
in $\pro{-}\MIC^\un(\dX,\cO_\dX)$, and for any second $K$-rational point~$y\in Y^\an(K)$ also a comparison map
\begin{equation}\label{eq:RH-et_functoriality}
	\pi_1^\RH(\dX;x,y) \to \sB_\deR\otimes_{\bQ_p}\pi_1^\ket(X^\an_{\bC_K};\bar x,\bar y)
\end{equation}
in the category of~$\sB_\deR$-schemes. These maps are compatible with all groupoid operations. Lemma~\ref{lem:etale-RH_extension_comparison} implies the following comparison result.

\begin{prop}\label{prop:etale-RH_comparison}
	The maps~\eqref{eq:RH-et_functoriality_universal} and~\eqref{eq:RH-et_functoriality} are isomorphisms for all points~$x,y\in Y^\an(K)$, and are~$G_K$-equivariant and strictly compatible with Hodge filtrations.
	\begin{proof}
		Using Lemma~\ref{lem:etale-RH_extension_comparison}, Proposition~\ref{prop:tannakian_iso} shows that the maps~\eqref{eq:RH-et_functoriality_universal} and~\eqref{eq:RH-et_functoriality} are isomorphisms. For $G_K$-equivariance, we note that~\eqref{eq:RH-et_functoriality_universal} takes~$e_x^\RH$ to~$1\otimes e^\ket_{\bar x}\in\dRH({}_{\bar x}\bar\bE^\ket)_x^{G_K} = (\sB_\deR\otimes_{\bQ_p}{}_{\bar x}\bar\bE^\ket)_x^{G_K}$. This implies by Lemma~\ref{lem:RH_equivariant_representability} that~\eqref{eq:RH-dR_functoriality_universal} is $G_K$-equivariant, and this automatically implies the same for~\eqref{eq:RH-et_functoriality}.
		
		As in the proof of Proposition~\ref{prop:dR-RH_comparison}, to show that~\eqref{eq:RH-et_functoriality_universal} and~\eqref{eq:RH-et_functoriality} are filtered isomorphisms, we show that the isomorphisms
		\[
		\phi_n\colon {}_x\dE^\RH_n\to\dRH({}_{\bar x}\bar\bE_n^\ket)
		\]
		of~\eqref{eq:tannakian_functoriality_family_levelwise} are filtered isomorphisms for all~$n$. The map~$\phi_n$ is certainly filtered by Proposition~\ref{prop:RH_filtered_representability}; we argue inductively that it is strict. This is obvious in the base case~$n=0$. Supposing then inductively that~$\phi_n$ is strict for the Hodge filtration, the proof of Proposition~\ref{prop:tannakian_iso} gives a commuting diagram
		\begin{center}
			\begin{tikzcd}
				0 \arrow[r] & \cO_\dX\otimes_{\sB_\deR}\rH^1_\deR(\dX,{}_x\dE^{\RH,\vee}_n)^\vee \arrow[r]\arrow[d,"\wr"'] & {}_x\dE^\RH_{n+1} \arrow[r]\arrow[d,"\phi_{n+1}","\wr"'] & {}_x\dE^\RH_n \arrow[r]\arrow[d,,"\phi_n","\wr"'] & 0 \\
				0 \arrow[r] & \cO_\dX\otimes_{\bQ_p}\rH^1_\ket(X^\an_{\bC_K},{}_{\bar x}\bar\bE^{\ket,\vee}_n)^\vee \arrow[r] & \dRH({}_{\bar x}\bar\bE^\ket_{n+1}) \arrow[r] & \dRH({}_{\bar x}\bar\bE^\ket_n) \arrow[r] & 0
			\end{tikzcd}
		\end{center}
		with strict exact rows, where the unmarked map is dual to the composite map
		\[
		\sB_\deR\otimes_{\bQ_p}\rH^1_\ket(X^\an_{\bC_K},{}_{\bar x}\bar\bE^{\ket,\vee}_n) \xrightarrow{\sim} \rH^1_\deR(\dX,\dRH({}_{\bar x}\bar\bE^{\ket,\vee}_n)) \xrightarrow[\sim]{\phi_n^*}  \rH^1_\deR(\dX,{}_x\dE^{\RH,\vee}_n) \,.
		\]
		By Lemma~\ref{lem:RH-dR_extension_comparison} both of the outer maps in the above diagram are strict for the filtrations. So by the five-lemma~$\phi_{n+1}$ is also strict, completing the induction.
	\end{proof}
\end{prop}

\subsection{Comparison of universal objects}\label{ss:comparison}

The Comparison Theorem~\ref{thm:comparison} is now easy to extract from Propositions~\ref{prop:dR-RH_comparison} and~\ref{prop:etale-RH_comparison}. For part~\eqref{thmpart:comparison_groupoid}, we have for every $K$-rational points~$x,y\in Y(K)$ a canonical isomorphism
\[
\sB_\deR\otimes_{\bQ_p}\pi_1^\ket(X^\an_{\bC_K};\bar x,\bar y) \xrightarrow\sim \sB_\deR\otimes_K\pi_1^\deR(X^\an;x,y)
\]
of affine $\sB_\deR$-schemes, $G_K$-equivariant and strictly compatible with the Hodge filtrations on affine rings. Moreover, these isomorphisms are compatible with the identities, composition and inversion in the fundamental groupoids. Thus, up to using different names for the objects involved, we have proved~\eqref{thmpart:comparison_groupoid}.

For the statement of~\eqref{thmpart:comparison_family}, we proceed somewhat indirectly. In our case, part~\eqref{thmpart:comparison_groupoid} implies in particular that the affine ring of~$\pi_1^\ket(X_{\bC_K};\bar x,\bar y)$ is an ind-de Rham representation of~$G_K$. So all of its finite-dimensional $G_K$-stable subrepresentations are de Rham. By Proposition~\ref{prop:tannakian_groupoid_via_universal_objects}, this implies that the fibre of~${}_x\bE^\et_n={}_x\bE^\ket_n|_{Y^\an_\et}$ at any $K$-rational point of~$Y^\an$ is de Rham for any~$n$, where~${}_x\bE^\ket_n$ is the descent of~${}_{\bar x}\bar\bE^\ket_n$ described in Lemma~\ref{lem:descent_of_universal}. Hence by \cite[Theorem~1.5(iii)]{liu-zhu:riemann-hilbert}, ${}_x\bE^\et_n$ is a de Rham local system on~$Y^\an_\et$ in the sense of Scholze~\cite[Definition~8.3]{scholze:relative}, for all~$n$.

This means that~${}_x\bE^\et_n$ has an associated filtered vector bundle with integrable connection on~$Y^\an$, which we denote by~$\sD_\deR({}_x\bE^\et_n)$. Recall that being \emph{associated} means that there exists an isomorphism
\[
c\colon \cO\bB_{\deR,Y^\an}\otimes_{\hat\bQ_p}{}_x\hat\bE^\et_n \xrightarrow\sim \cO\bB_{\deR,Y^\an}\otimes_{\mu^{-1}\cO_{Y^\an}}\mu^{-1}\sD_\deR({}_x\bE^\et_n)
\]
of filtered $\cO\bB_{\deR,Y^\an}$-vector bundles with integrable connection on~$Y^\an_\proet$, where~$\mu\colon Y^\an_\proet\to Y^\an_\an$ is the projection to the analytic site\footnote{Actually, Scholze uses the \'etale site of~$Y^\an$ instead of the analytic site, but this makes no difference due to \cite[Lemma~7.3]{scholze:relative}.}. It follows from \cite[Theorem~7.6(i) \& Lemma~7.3]{scholze:relative} that $\sD_\deR({}_x\bE^\et_n)$ has the following canonical description: we have
\[
\sD_\deR({}_x\bE^\et_n) = \mu_*(\cO\bB_{\deR,Y^\an}\otimes_{\hat\bQ_p}\hat\bE)
\]
and the isomorphism~$c$ is the inverse of the map induced by the pull--push adjunction.

On the other hand, we can also consider
\[
\sD_\deR^\wedge({}_x\bE^\ket_n) = \mu_*(\cO\bB_{\deR,X^\an}^\wedge\otimes_{\hat\bQ_p}\hat\bE) \,,
\]
where~$\cO\bB_{\deR,X^\an}^\wedge$ is the period sheaf of \cite[Definition~2.2.10]{dllz:riemann-hilbert} (which we were denoting simply~$\cO\bB_\deR^\wedge$ earlier) and~$\mu\colon X^\an_\proket\to X^\an_\an$ is the natural morphism of sites (by mild abuse of notation). This object~$\sD_\deR^\wedge({}_x\bE^\ket_n)$ is a vector bundle with integrable logarithmic connection on~$X^\an$, and is endowed with a Griffiths-transverse filtration by coherent subsheaves \cite[Theorem~1.7]{dllz:riemann-hilbert}. Since~$\cO\bB_{\deR,X^\an}^\wedge|_{Y^\an_\proet}$ is the completion of~$\cO\bB_{\deR,Y^\an}$ with respect to its filtration (see \cite[Remark~2.2.11]{dllz:riemann-hilbert}), we see that~$\sD_\deR^\wedge({}_x\bE^\ket_n)|_{Y^\an_\an}$ is the completion of~$\sD_\deR({}_x\bE^\et_n)$ with respect to its filtration, so we have
\[
\sD_\deR({}_x\bE^\et_n) = \sD_\deR^\wedge({}_x\bE^\ket_n)|_{Y^\an_\an} \,.
\]

Now it follows from the proof of \cite[Lemma~3.3.17]{dllz:riemann-hilbert} that we have
\[
\sD_\deR^\wedge({}_x\bE^\ket_n) = \dRH({}_{\bar x}\bar\bE^\ket_n)^{G_K} \,,
\]
while Propositions~\ref{prop:dR-RH_comparison} and~\ref{prop:etale-RH_comparison} provide us with $G_K$-equivariant filtered isomorphisms
\begin{equation}\label{eq:composite_comparison}
	\dRH({}_{\bar x}\bar\bE^\ket_n) \cong {}_x\dE^\RH_n \cong \sB_\deR\hatotimes_K{}_x\dE^\deR_n \,.
\end{equation}
Combining all of these, we have constructed an isomorphism
\[
\sD_\deR({}_x\bE^\et_n)\cong{}_x\dE^\deR_n|_{Y^\an_\an}
\]
for all~$n$. It is easy to check (using $\otimes$-naturality of the constructions involved) that these isomorphisms are compatible as~$n$ varies, so induce a canonical isomorphism
\[
\sD_\deR({}_x\bE^\et)\cong{}_x\dE^\deR|_{Y^\an_\an}
\]
of coalgebras in $\pro{-}\MIC^\un(Y^\an,\cO_{Y^\an})$. So~${}_x\bE^\et$ and~${}_x\dE^\deR|_{Y^\an_\an}$ are associated in a canonical manner, meaning that we have the isomorphism~\eqref{eq:comparison_family} claimed in part~\eqref{thmpart:comparison_family}.\qed

\begin{rmk}
	One can also run a direct proof, starting from the isomorphism~\eqref{eq:composite_comparison} and transforming it into the desired isomorphism~\eqref{eq:comparison_family}; deducing \emph{a posteriori} that~${}_x\bE^\et_n$ is de Rham. However, due to the differences in conventions between~\cite{scholze:relative} and~\cite{dllz:riemann-hilbert} (chiefly that the period sheaves in~\cite{dllz:riemann-hilbert} are the completions of their counterparts in~\cite{scholze:relative}), this direct argument is rather fiddly to write down carefully.
\end{rmk}
\section{Horizontal de Rham local systems over polydiscs}

The other ingredient we will need in order to prove Theorem~\ref{thm:local_constancy_of_Kummer_maps_p-adic} is a theorem of Shimizu on potential horizontal semistability for horizontal de Rham local systems over spherical polyannuli. The exposition here roughly follows \cite[\S3.5]{me-jakob:lawrence-venkatesh_sections}. Let~$V$ be a rigid-analytic space over~$K$ isomorphic to either a closed polydisc $\Sp(K\langle t_1,\dots,t_n\rangle)$ or a spherical polyannulus $\Sp(K\langle t_1^{\pm1},\dots,t_n^{\pm1}\rangle)$. A vector bundle~$\dE$ with integrable connection~$\nabla$ on~$V$ is said to have a \emph{full basis of horizontal sections} just when there is an isomorphism
\[
\phi\colon (\cO_V,\rd)^{\oplus m} \xrightarrow\sim (\dE,\nabla)
\]
for some~$m$. Since~$\rH^0_\deR(V,\cO_V)=K$, it follows that the choice of isomorphism~$\phi$ is unique up to the action of~$\GL_m(K)$. Hence for any two points~$x,y\in V(K)$ there is a well-defined $K$-linear \emph{parallel transport} isomorphism
\[
T_{x,y}^\nabla\colonequals \phi_y\circ\phi_x^{-1} \colon \dE_x \xrightarrow\sim \dE_y \,.
\]


The precise result we need to extract from Shimizu's theory is the following.

\begin{thm}\label{thm:parallel_transport_is_phi-compatible}
	Let~$V$ be a rigid-analytic space over~$K$ isomorphic to a closed polydisc or spherical polyannulus. Let~$\bE$ be an isogeny $\bZ_p$-local system on~$V_\et$, $\dE$ a filtered vector bundle with integrable connection on~$V$, and
	\[
	c\colon \cO\bB_{\deR,V}\otimes_{\hat\bQ_p}\hat\bE \xrightarrow\sim \cO\bB_{\deR,V}\otimes_{\mu^{-1}\cO_V}\mu^{-1}\dE
	\]
	an isomorphism of filtered $\cO\bB_{\deR,V}$-vector bundles with integrable connection on~$V_\proet$ (so~$\bE$ is de Rham, associated to~$\dE$). Suppose that~$\dE$ has a full basis of horizontal sections. Then for every $x,y\in V(K)$, there is a unique isomorphism
	\[
	T_{x,y}\colon\sD_\pst\left(\bE_{\bar x}\right) \xrightarrow\sim \sD_\pst\left(\bE_{\bar y}\right)
	\]
	of $(\varphi,N,G_K)$-modules making the diagram
	\begin{equation}\label{diag:parallel_transport_is_phi-compatible}
	\begin{tikzcd}[column sep = large]
		\sB_\deR\otimes_{\bQ_p^\nr}\sD_\pst\left(\bE_{\bar x}\right) \arrow[r,"1\otimes\alpha_\pst","\sim"']\arrow[d,"1\otimes T_{x,y}"',"\wr"] & \sB_\deR\otimes_{\bQ_p}\bE_{\bar x} \arrow[r,"1\otimes c_x","\sim"'] & \sB_\deR\otimes_K\dE_x \arrow[d,"1\otimes T_{x,y}^\nabla","\wr"'] \\
		\sB_\deR\otimes_{\bQ_p^\nr}\sD_\pst\left(\bE_{\bar y}\right) \arrow[r,"1\otimes\alpha_\pst","\sim"'] & \sB_\deR\otimes_{\bQ_p}\bE_{\bar y} \arrow[r,"1\otimes c_y","\sim"'] & \sB_\deR\otimes_K\dE_y
	\end{tikzcd}
	\end{equation}
	commute. Here $\sD_\pst$ denotes Fontaine's potentially semistable functor \cite[\S5.6.4]{fontaine3}, and $\alpha_\pst$ denotes the natural map
	\[
	\sB_\st\otimes_{\bQ_p^\nr}\sD_\pst(E) \to \sB_\st\otimes_{\bQ_p}E
	\]
	for any $\bQ_p$-linear representation~$E$ of~$G_K$ (which is an isomorphism for~$E$ de Rham).
\end{thm}

\begin{rmk}
	In order to make sense of the arrows marked~$1\otimes\alpha_\pst$ in the above diagram, one has to choose an embedding~$\sB_\st\hookrightarrow\sB_\deR$, i.e.\ a branch of the $p$-adic logarithm. We fix such a choice once and for all. Although we will not need this, we remark that the map~$T_{x,y}$ of Theorem~\ref{thm:parallel_transport_is_phi-compatible} is independent of this choice.
\end{rmk}

\begin{rmk}\label{rmk:parallel_transport_tensor-natural}
	The isomorphism~$T_{x,y}$ appearing in Theorem~\ref{thm:parallel_transport_is_phi-compatible} is $\otimes$-natural in~$\bE$, since the parallel transport map~$T_{x,y}^\nabla$ is $\otimes$-natural in~$\dE$.
\end{rmk}

We first prove Theorem~\ref{thm:parallel_transport_is_phi-compatible} in the case that~$V=\bT^n:=\Sp(R[p^{-1}])$ is the spherical polyannulus with $R:=\cO_K\langle t_1^{\pm1},t_2^{\pm1},\dots,t_n^{\pm1}\rangle$. Fix a geometric generic point~$\bar\eta$ of~$V$, by which we mean a $K$-embedding $\bar\eta\colon R\hookrightarrow\Omega$ into an algebraically closed field~$\Omega$ containing~$\Kbar$. We write~$\Rbar\subseteq\Omega$ for the union of the finite $R$-subalgebras $R'\subseteq\Omega$ such that~$R'[p^{-1}]$ is \'etale over~$R[p^{-1}]$. So~$\Rbar[p^{-1}]$ is the union of all finite \'etale $R[p^{-1}]$-subalgebras of~$\Omega$, and~$\Rbar$ is the normalisation of~$R$ in the fraction field of~$\Rbar[p^{-1}]$. We write
\[
G_{R_K}:=\Aut(\Rbar[p^{-1}]/R[p^{-1}]) \,.
\]
By construction, we have $\Kbar\subseteq\Rbar[p^{-1}]$, and restriction of the tautological action provides a surjective homomorphism
\[
G_{R_K}\twoheadrightarrow G_K \,.
\]
If~$L$ is a finite extension of~$K$ contained in~$\Kbar$, we write~$G_{R_L}\subseteq G_{R_K}$ for the preimage of~$G_L$ under the restriction map.

The group~$G_{R_K}$ is canonically identified with the profinite \'etale fundamental group of $\Spec(R[p^{-1}])$ based at~$\bar\eta$ \cite[Proposition~5.4.9]{szamuely:fundamental_groups}. This is the same as the algebraic fundamental group $\pi_1^\alg(V,\bar\eta)$ of~$V$ in the sense of de Jong \cite[p.~94]{de_jong:rigid_pi_1}, as follows.

\begin{lem}\label{lem:analytic_and_algebraic_pi_1}
	Let~$A$ be a strictly affinoid $K$-algebra, with associated affinoid $K$-analytic space~$\Sp(A)$. Then the functor $B\mapsto\Sp(B)$ provides a contravariant equivalence of categories between finite \'etale $A$-algebras and finite \'etale coverings of~$\Sp(A)$.
	
	In particular, if~$A$ is an integral domain and $\bar\eta\colon A\to\Omega$ is a homomorphism from~$A$ into a complete algebraically closed field, then there is a canonical isomorphism
	\[
	\pi_1^\alg(\Sp(A),\bar\eta) \cong \pi_1^\et(\Spec(A),\bar\eta) \,.
	\]
	\begin{proof}
		\cite[Definition~4.5.7]{fresnel-van_der_put:rigid} and the result just below it show that $B\mapsto\Sp(B)$ gives an equivalence between finite $A$-algebras and $K$-analytic spaces finite over~$\Sp(A)$. Part~(1) of the proof of \cite[Proposition~8.1.1]{fresnel-van_der_put:rigid} shows that a morphism $\Sp(B)\to\Sp(A)$ of affinoid $K$-analytic spaces is \'etale if and only if~$B$ is a flat $A$-algebra and the module $\Omega_{B/A}^{1,f}$ of finite differentials is zero. In the particular case that~$B$ is a finite $A$-algebra, the module~$\Omega_{B/A}^1$ of K\"ahler differentials is finite, whence $\Omega_{B/A}^{1,f}=\Omega_{B/A}^1$. We conclude that a finite morphism $\Sp(B)\to\Sp(A)$ is \'etale just when~$B$ is flat over~$A$ and $\Omega^1_{B/A}=0$; that is to say, just when~$B$ is an \'etale $A$-algebra.
	\end{proof}
\end{lem}

Now associated to~$\Rbar[p^{-1}]$ is a certain object of the pro-\'etale site of~$V$, namely
\[
\Vbar:=\varprojlim_{R'}\Sp(R'[p^{-1}]) \,,
\]
where $R'$ runs over all finite $R$-algebras in~$\Omega$ such that $R'[p^{-1}]$ is \'etale over~$R[p^{-1}]$. There is a tautological action of~$G_{R_K}$ on~$\Vbar$ from the right, which is the same as the usual action of the fundamental group~$G_{R_K}$ on the pointed universal covering~$\Vbar$. Following \cite[p.~47]{shimizu:p-adic_monodromy}, we write\footnote{Contrary to what the notation suggests here, $\sB_\deR(R)$ depends on~$\Rbar$, not just~$R$.}
\[
\sB_\deR(R):=\cO\bB_{\deR,V}(\Vbar) \,,
\]
which is a filtered $\Rbar[p^{-1}]$-algebra endowed with an action of~$G_{R_K}$ extending the tautological action on~$\Rbar[p^{-1}]$, and with a connection
\[
\nabla\colon \sB_\deR(R) \to \Omega^{1,f}_{R[p^{-1}]/K}\otimes_{R[p^{-1}]}\sB_\deR(R)
\]
satisfying the Leibniz rule with respect to the derivation on~$R[p^{-1}]$. We write~$\sB_\deR^\nabla(R)$ for the kernel of the connection. The $G_{R_K}$-invariant subring of~$\sB_\deR^\nabla(R)$ is~$K$ \cite[Proposition~4.9]{shimizu:p-adic_monodromy}.

Shimizu defines a $G_{R_K}$-stable subring $\sB_\st^\nabla(R)\subseteq\sB_\deR^\nabla(R)$ possessing a semilinear crystalline Frobenius endomorphism~$\varphi$ and monodromy operator~$N$ satisfying all the usual relations \cite[Definition~4.7]{shimizu:p-adic_monodromy}. The $G_{R_K}$-invariant subring of  $\sB_\st^\nabla$ is the maximal absolutely unramified subextension~$K_0$ of~$K$ \cite[Corollary~4.10]{shimizu:p-adic_monodromy}; more generally, we have $(\sB_\st^\nabla)^{G_{R_L}}=L_0$ for every finite extension~$L$ of~$K$ contained in~$\Kbar$.

In this setup, we have analogues of Fontaine's functors for representations of~$G_{R_K}$.

\begin{defi}[{\cite[Definitions~4.7 \&~4.16]{shimizu:p-adic_monodromy}}]
	Let~$E$ be a continuous representation of~$G_{R_K}$ on a $\bQ_p$-vector space. Then we define
	\[
	\sD_\deR^\nabla(E) \colonequals (\sB_\deR^\nabla(R)\otimes_{\bQ_p}E)^{G_{R_K}} \hspace{0.4cm}\text{and}\hspace{0.4cm} \sD_\pst^\nabla(E)\colonequals\varinjlim_L(\sB_\st^\nabla(R)\otimes_{\bZ_p}E)^{G_{R_L}}
	\]
	where the colimit is taken over finite extensions $L/K$ contained in~$\Kbar$. The natural action of $G_{R_K}$ on~$\sD_\pst^\nabla(E)$ factors through~$G_K$ and has finite point-stabilisers: together with the crystalline Frobenius~$\varphi$ and monodromy operator~$N$ induced from those on~$\sB_\st^\nabla(R)$, this makes $\sD_\pst^\nabla(E)$ into a discrete $(\varphi,N,G_K)$-module.
	
	We always have the inequality $\dim_K\sD_\deR(E)\leq\dim_{\bQ_p}E$ \cite[Lemma~4.12]{shimizu:p-adic_monodromy}, and we say that~$E$ is \emph{horizontal de Rham} just when equality holds. Equivalently, $E$ is horizontal de Rham just when the evident map
	\[
	\alpha_\deR^\nabla\colon \sB_\deR^\nabla(R)\otimes_K\sD_\deR^\nabla(E) \to \sB_\deR^\nabla(R)\otimes_{\bQ_p}E
	\]
	is a $\sB_\deR^\nabla(R)$-linear isomorphism. (One can also define potentially horizontal semistable representations via a similar formalism, but it turns out that every horizontal de Rham representation is potentially horizontal semistable and vice versa, as we will see shortly.)
\end{defi}

This theory is contravariant functorial in the ring~$\Rbar$ \cite[\S4.3]{shimizu:p-adic_monodromy}. Suppose that~$x\in V(K)$ is a $K$-point of~$V$, and choose an extension of the map $x^*\colon R[p^{-1}]\to K$ to a $\Kbar$-algebra homomorphism $\tilde x^*\colon\Rbar[p^{-1}] \to \Kbar$. There are then induced maps
\[
\tilde x^*\colon \sB_\st^\nabla(R)\to\sB_\st^\nabla(\cO_K)=\sB_\st \hspace{0.4cm}\text{and}\hspace{0.4cm} \tilde x^*\colon \sB_\deR(R)\to\sB_\deR(\cO_K)=\sB_\deR \,,
\]
the latter of which is the map $\rH^0(\Vbar,\cO\bB_{\deR,V})\to\rH^0(\Sp(\Kbar),\cO\bB_{\deR,\Sp(K)})$ given by taking the fibre at~$x$. These maps are compatible with all of the usual structures and inclusions between these algebras.

These maps~$\tilde x^*$ are also Galois-equivariant, in the following sense. The $\Kbar$-algebra homomorphism $\tilde x^*$ determines a $\Kbar$-point~$\tilde x$ on~$\Vbar$ lying over $x\in V(K)$. This determines in turn a profinite \'etale path $\gamma_{\tilde x}\in\pi_1^\alg(V;\bar\eta,\bar x)$, and hence a splitting~$s_{\tilde x}$ of the map $G_{R_K}\to G_K$. This splitting is characterised by the fact that map $\tilde x^*\colon\Rbar[p^{-1}]\to\Kbar$ is $G_K$-equivariant for the action on~$\Rbar[p^{-1}]$ given by restricting the $G_{R_K}$-action along~$s_{\tilde x}$. Since the construction of $\sB_\st^\nabla(R)$ and $\sB_\deR(R)$ is functorial in~$\Rbar[p^{-1}]$, this also implies that the maps $\tilde x^*\colon\sB_\st^\nabla(R)\to\sB_\st$ and $\tilde x^*\colon\sB_\deR(R)\to\sB_\deR$ are both $G_K$-equivariant in the same sense.

If now~$E$ is a continuous representation of~$G_{R_K}$ on a $\bQ_p$-vector space, we write~$E|_{\tilde x}$ for~$E$ endowed with the $G_K$-action given by restriction along~$s_{\tilde x}$. There are natural maps
\[
\psi^\deR_{\tilde x}\colon\sD_\deR^\nabla(E) \to \sD_\deR(E|_{\tilde x}) \hspace{0.4cm}\text{and}\hspace{0.4cm} \psi^\pst_{\tilde x}\colon\sD_\pst^\nabla(E) \to \sD_\pst(E|_{\tilde x})
\]
given by taking invariants in $\tilde x^*\otimes1\colon\sB_\deR^\nabla(R)\otimes_{\bZ_p}E \to \sB_\deR\otimes_{\bQ_p}E|_{\tilde x}$ and $\tilde x^*\otimes1\colon\sB_\st^\nabla(R)\otimes_{\bZ_p}E \to \sB_\st\otimes_{\bQ_p}E|_{\tilde x}$. The map~$\psi^\deR_{\tilde x}$ is $K$-linear, and~$\psi^\pst_{\tilde x}$ is a morphism of discrete~$(\varphi,N,G_K)$-modules. We summarise the key results from Shimizu's paper which we will need.

\begin{prop}\label{prop:shimizu_results}
	Suppose that~$E$ is a horizontal de Rham representation of~$G_{R_K}$. Then:
	\begin{enumerate}
		\item $\dim_{\bQ_p^\nr}\sD_\pst^\nabla(E) = \dim_{\bQ_p}E$ and the natural map
		\[
		\Kbar\otimes_{\bQ_p^\nr}\sD_\pst^\nabla(E) \to \Kbar\otimes_K\sD_\deR^\nabla(E)
		\]
		is an isomorphism;
		\item for any point~$x\in V(K)$ and lift~$\tilde x\in\Vbar(\Kbar)$, the restricted representation~$E|_{\tilde x}$ is de Rham, and the maps~$\psi^\deR_{\tilde x}$ and~$\psi^\pst_{\tilde x}$ are isomorphisms fitting into a commuting square
		\begin{center}
		\begin{tikzcd}
			\Kbar\otimes_{\bQ_p^\nr}\sD_\pst^\nabla(E) \arrow[r,"\sim"]\arrow[d,"1\otimes\psi^\pst_{\tilde x}","\wr"'] & \Kbar\otimes_K\sD_\deR^\nabla(E) \arrow[d,"1\otimes\psi^\deR_{\tilde x}","\wr"'] \\
			\Kbar\otimes_{\bQ_p^\nr}\sD_\pst(E|_{\tilde x}) \arrow[r,"\sim"] & \Kbar\otimes_K\sD_\deR(E|_{\tilde x}) \,.
		\end{tikzcd}
		\end{center}
	\end{enumerate}
	\begin{proof}
		For the first point, Shimizu's $p$-adic monodromy theorem \cite[Theorem~7.4]{shimizu:p-adic_monodromy} implies that~$E$ is potentially horizontal semistable, meaning that the natural map $\sB_\st^\nabla(R)\otimes_{\bQ_p^\nr}\sD_\pst^\nabla(E) \to \sB_\st^\nabla(R)\otimes_{\bQ_p}E$ is an isomorphism. This directly gives the equality of dimensions, and by tensoring up to~$\sB_\deR^\nabla(R)$ and taking $G_{R_L}$-fixed points, also the isomorphy of the claimed map. In the second point, the fact that~$E|_{\tilde x}$ is de Rham and that~$\psi^\deR_{\tilde x}$ is an isomorphism is \cite[Lemma~4.4]{shimizu:p-adic_monodromy}. That the square commutes is obvious, and this then implies that~$\psi^\pst_{\tilde x}$ is also an isomorphism.
	\end{proof}
\end{prop}

In particular, if~$y\in V(K)$ is a second $K$-rational point and~$\tilde y$ is a lift, then there are natural isomorphisms
\[
T_{\tilde x,\tilde y}^\deR = \psi^\deR_{\tilde y}\circ\psi^{\deR,-1}_{\tilde x} \colon \sD_\deR(E|_{\tilde x}) \xrightarrow\sim \sD_\deR(E|_{\tilde y}) \hspace{0.4cm}\text{and}\hspace{0.4cm} T_{\tilde x,\tilde y}^\pst = \psi^\pst_{\tilde y}\circ\psi^{\pst,-1}_{\tilde x} \colon \sD_\pst(E|_{\tilde x}) \xrightarrow\sim \sD_\pst(E|_{\tilde y})
\]
for any horizontal de Rham representation~$E$, the latter of which is an isomorphism of discrete~$(\varphi,N,G_K)$-modules.

\smallskip

Using this, we are in a position to prove Theorem~\ref{thm:parallel_transport_is_phi-compatible}, beginning with the case that~$V$ is a spherical polyannulus. We set~$E=\rH^0(\Vbar,\bE)=\bE_{\bar\eta}$ and~$\cE=\rH^0(V,\dE)$, so that we have $E|_{\tilde x}=\bE_{\bar x}$ and~$E|_{\tilde y}=\bE_{\bar y}$, and~$\dE_x$ and~$\dE_y$ are the base-changes of~$\cE$ along~$x^*,y^*\colon R[p^{-1}]\to K$, respectively. We will show that the map~$T^\pst_{\tilde x,\tilde y}$ above makes the desired square commute; it is certainly then unique.

Taking sections over~$\Vbar$ of the isomorphism $c\colon\cO\bB_{\deR,V}\otimes_{\hatotimes\bQ_p}\hat\bE\xrightarrow\sim\cO\bB_{\deR,V}\otimes_{\cO_V}\dE$ yields a $G_{R_K}$-equivariant isomorphism
\begin{equation}\label{eq:sections_of_comparison}
\sB_\deR(R)\otimes_{\bQ_p}E\xrightarrow\sim\sB_\deR(R)\otimes_{R[p^{-1}]}\cE
\end{equation}
of~$\sB_\deR(R)$-modules, compatible with connections. The right-hand side is isomorphic by assumption to $\sB_\deR(R)^{\oplus m}$ for some~$m$. Since~$\sB_\deR^\nabla(E)^{G_{R_K}}=K$ by \cite[Proposition~4.9]{shimizu:p-adic_monodromy}, we obtain from~\eqref{eq:sections_of_comparison} an isomorphism
\[
\sD_\deR(E)\xrightarrow\sim\cE^{\nabla=0} \,.
\]
In particular, $\dim_K\sD_\deR(E)=\dim_K(\cE^{\nabla=0})=\dim_{\bQ_p}E$ and so~$E$ is a horizontal de Rham representation.

Moreover, the induced map
\[
\phi\colon R[p^{-1}]\otimes_K\sD_\deR^\nabla(E) \to \cE
\]
is an isomorphism, giving a trivialisation of~$\cE$ as an $R[p^{-1}]$-module with connection. Using this trivialisation to calculate the parallel transport, we obtain a commuting diagram
\begin{center}
\begin{tikzcd}
	\sD_\deR(E|_{\tilde x}) \arrow[d,"T_{\tilde x,\tilde y}^\deR","\wr"'] & \sD_\deR^\nabla(E) \arrow[d,equals]\arrow[l,"\psi_{\tilde x}^\deR"',"\sim"]\arrow[r,"\phi_x","\sim"'] & \dE_x \arrow[d,"T_{x,y}^\nabla","\wr"'] \\
	\sD_\deR(E|_{\tilde y}) & \sD_\deR^\nabla(E) \arrow[l,"\psi_{\tilde y}^\deR"',"\sim"]\arrow[r,"\phi_y","\sim"'] & \dE_y
\end{tikzcd}
\end{center}
Tensoring up with~$\sB_\deR$ then gives commutativity of
\begin{center}
\begin{tikzcd}
	\sB_\deR\otimes_K\sD_\deR(E|_{\tilde x}) \arrow[d,"1\otimes T_{\tilde x,\tilde y}^\deR","\wr"']\arrow[r,"\alpha_\deR","\sim"'] & \sB_\deR\otimes_{\bQ_p}E|_{\tilde x} \arrow[r,"c_x","\sim"'] & \sB_\deR\otimes_K\dE_x \arrow[d,"1\otimes T_{x,y}^\nabla","\wr"'] \\
	\sB_\deR\otimes_K\sD_\deR(E|_{\tilde y}) \arrow[r,"\alpha_\deR","\sim"'] & \sB_\deR\otimes_{\bQ_p}E|_{\tilde y} \arrow[r,"c_y","\sim"'] & \sB_\deR\otimes_K\dE_y
\end{tikzcd}
\end{center}
using that the fibre~$c_x$ of the isomorphism~$c$ is the base-change of~\eqref{eq:sections_of_comparison} along the map $\tilde x^*\colon\sB_\deR(R)\to\sB_\deR$. Combined with Proposition~\ref{prop:shimizu_results}, this implies that the map~$T_{x,y}^\pst=T_{\tilde x,\tilde y}^\pst$ makes~\eqref{diag:parallel_transport_is_phi-compatible} commute. It is clearly the unique such map, and so we are done.

It remains to deal with the case that~$V$ is a polydisc, for which we reduce to the case of a spherical polyannulus. When the residue field of~$K$ is not~$\bF_2$, we may choose a spherical polyannulus~$V^\circ\subseteq V$ containing~$x$ and~$y$, and the result for~$V$ follows from the result for~$V^\circ$. If the residue field is~$\bF_2$, we let~$K_d$ for~$d>1$ denote the unramified field extension of degree~$d$ inside~$\Kbar$. The result for~$V_{K_d}$ implies that there is a unique isomorphism~$T_{x,y}^{\pst,d}\colon\sD_\pst(\bE_{\bar x})\xrightarrow\sim\sD_\pst(\bE_{\bar y})$ of~$(\varphi,N,G_{K_d})$-modules making~\eqref{diag:parallel_transport_is_phi-compatible} commute. Unicity implies that in fact all the~$T_{x,y}^{\pst,d}$ are equal, so define a $\bQ_p^\nr$-linear isomorphism $T_{x,y}^\pst\colon\sD_\pst(\bE_{\bar x})\xrightarrow\sim\sD_\pst(\bE_{\bar y})$ which is compatible with the Frobenius~$\varphi$, monodromy operator~$N$, and the action of~$G_{K_d}\subseteq G_K$ for any~$d>1$. Since~$G_{K_2}$ and~$G_{K_3}$ generate~$G_K$, we see that~$T_{x,y}^\pst$ is $G_K$-equivariant and we are done. \qed
\section{Local constancy of the $\bQ_p$-pro-unipotent Kummer map mod~$\rH^1_e$}
\label{s:main_theorem_p-adic}

Now we are in a position to prove Theorem~\ref{thm:local_constancy_of_Kummer_maps_p-adic}. Before we do, we need a single input.

\begin{prop}\label{prop:enough_flat_sections}
	Let~$V\subseteq Y^\an$ be an admissible open, isomorphic to a closed polydisc. Then~$\dE|_V$ has a full basis of horizontal sections for every~$\dE\in\MIC^\un(Y^\an,\cO_{Y^\an})$.
	\begin{proof}
		If we identify~$V$ with the closed polydisc of radii~$1$, then there is some~$\rho>1$ such that the inclusion~$V\hookrightarrow Y^\an$ extends to an inclusion~$V_\rho\hookrightarrow Y^\an$ of a closed polydisc of radii~$\rho$, necessarily uniquely. We let
		\[
		\MIC^\un(V,\cO_V^\dagger) = \varinjlim_{\rho>1}\MIC^\un(V_\rho,\cO_{V_\rho})
		\]
		denote the category of overconvergent vector bundles with integrable connection on~$V$, so that the restriction map $\MIC^\un(Y^\an,\cO_{Y^\an})\to\MIC^\un(V,\cO_V)$ factors through $\MIC^\un(V,\cO_V^\dagger)$. Extensions in~$\MIC^\un(V,\cO_V^\dagger)$ are computed by overconvergent de Rham cohomology
		\[
		\Ext^1_{\MIC^\dagger}(\cO_V^\dagger,\dE^\dagger) = \rH^1_\deR(V,\dE^\dagger) \colonequals \bH^1(V,\DeR(\dE^\dagger))
		\]
		where $\DeR(\dE^\dagger)$ is the overconvergent de Rham complex, cf. \cite[Definition~7.5.11]{fresnel-van_der_put:rigid}. But $\rH^1_\deR(V,\cO_V^\dagger)=0$, so all objects in~$\MIC^\un(V,\cO_V^\dagger)$ are direct sums of the unit object, and the result follows.
	\end{proof}
\end{prop}

\begin{proof}[Proof of Theorem~\ref{thm:local_constancy_of_Kummer_maps_p-adic}]
	We will prove that if~$V\subseteq Y^\an$ is admissible open, isomorphic to a closed polydisc, then the unipotent Kummer map modulo~$\rH^1_e$ is constant on~$V$. Let~$x,y\in V(K)$, and write~${}_{\bar x}\bar\bE^\et$ for the universal object of~$\Loc^\un(Y_{\bC_K},\bQ_p)$ based at~$\bar x$. Let~${}_x\bE^\et$ be the descent of~${}_{\bar x}\bar\bE^\et$ to a local system on~$Y_\et$ as in Lemma~\ref{lem:descent_of_universal}; this is a pro-de Rham local system with associated vector bundle with connection~${}_x\dE^\deR$ by Theorem~\ref{thm:comparison}. Proposition~\ref{prop:enough_flat_sections} implies that~${}_x\dE^\deR|_V$ has a full basis of horizontal sections, so we deduce from Theorem~\ref{thm:parallel_transport_is_phi-compatible} an isomorphism
	\[
	T_{x,y}\colon \sD_\pst({}_{\bar x}\bar\bE^\et_{\bar x}) \xrightarrow\sim \sD_\pst({}_{\bar x}\bar\bE^\et_{\bar y})
	\]
	of pro-$(\varphi,N,G_K)$-modules. This isomorphism is compatible with the coalgebra structure on either side by Remark~\ref{rmk:parallel_transport_tensor-natural}.
	
	Via Proposition~\ref{prop:tannakian_groupoid_via_universal_objects}, we obtain a $G_K$-equivariant isomorphism
	\[
	\sB_\st\otimes_{\bQ_p^\nr}\cO(\pi_1^\et(Y_{\bC_K};\bar x,\bar y)) \xrightarrow\sim \sB_\st\otimes_{\bQ_p^\nr}\cO(\pi_1^\et(Y_{\bC_K};\bar x,\bar x))
	\]
	of~$\sB_\st$-algebras, compatible with the Frobenius automorphism and monodromy operator on either side. In particular, passing to $(\varphi,N)$-fixed points yields a $G_K$-equivariant isomorphism
	\[
	\sB_\cris^{\varphi=1}\otimes_{\bQ_p^\nr}\cO(\pi_1^\et(Y_{\bC_K};\bar x,\bar y)) \xrightarrow\sim \sB_\cris^{\varphi=1}\otimes_{\bQ_p^\nr}\cO(\pi_1^\et(Y_{\bC_K};\bar x,\bar x))
	\]
	of~$\sB_\cris^{\varphi=1}$-algebras. Hence
	\[
	\pi_1^\et(Y_{\bC_K};\bar x,\bar y)(\sB_\cris^{\varphi=1})^{G_K} \neq \emptyset \,,
	\]
	since it contains an element corresponding to~$1\in 	\pi_1^\et(Y_{\bC_K};\bar x,\bar x)(\sB_\cris^{\varphi=1})^{G_K}$. Composition with this element gives a $G_K$-equivariant isomorphism
	\[
	\pi_1^\et(Y_{\bC_K};\bar x_0,\bar x)(\sB_\cris^{\varphi=1}) \cong \pi_1^\et(Y_{\bC_K};\bar x_0,\bar y)(\sB_\cris^{\varphi=1})
	\]
	of $\pi_1^\et(Y_{\bC_K},\bar x_0)(\sB_\cris^{\varphi=1})$-torsors, whence~$j(x)=j(y)$ modulo~$\rH^1_e$ as desired.
\end{proof}

\bibliography{references}

\begin{thebibliography}{BDCKW14}

\bibitem[AGV72]{sga4}
M.~Artin, A.~Grothendieck, and J.L. Verdier.
\newblock {\em SGA 4: Th{\'e}orie des {T}opos et {C}ohomologie {\'E}tale des
  {S}ch{\'e}mas}, volume 269 of {\em Lecture Notes in Mathematics}.
\newblock Springer Verlag, 1972.

\bibitem[AIK15]{andreatta-iovita-kim}
F.~Andreatta, A.~Iovita, and M.~Kim.
\newblock A $p$-adic non-abelian criterion for good reduction of curves.
\newblock {\em Duke Mathematical Journal}, 164(13):2597--2642, 2015.

\bibitem[BDCKW14]{minhyong-etal:bsd_conjecture}
J.S. Balakrishnan, I.~Dan-Cohen, M.~Kim, and S.~Wewers.
\newblock A non-abelian conjecture of {Birch} and {Swinnerton-Dyer} type for
  hyperbolic curves, 2014.
\newblock arXiv: 1209.0640 (preprint only) -- v3 accessed.

\bibitem[Ber93]{berkovich}
V.~Berkovich.
\newblock {\'E}tale cohomology for non-archimedean analytic spaces.
\newblock {\em Publications math{\'e}matiques de l'I.H.{\'E}.S.}, 78:5--161,
  1993.

\bibitem[Bet18]{thesis}
L.A. Betts.
\newblock {\em Heights via anabelian geometry and local Bloch--Kato Selmer
  sets}.
\newblock PhD thesis, University of Oxford, 2018.

\bibitem[Bet19]{me:motivic_heights}
L.A. Betts.
\newblock The motivic anabelian geometry of local heights on abelian varieties,
  2019.
\newblock arXiv preprint (arXiv:1706.04850v2).

\bibitem[BL]{me-daniel:weight-monodromy}
L.A. Betts and D.~Litt.
\newblock Semisimplicity of the {F}robenius action on $\pi_1$.
\newblock arXiv:1912.02167 (v2 accessed).

\bibitem[BS22]{me-jakob:lawrence-venkatesh_sections}
L.A. Betts and J.~Stix.
\newblock Galois sections and $p$-adic period mappings, 2022.
\newblock arXiv:2204.13674 -- v1 accessed.

\bibitem[Del97]{deligne:local_behaviour}
P.~Deligne.
\newblock Local behavior of {H}odge structures at infinity.
\newblock {\em {AMS/IP} {S}tudies in {A}dvanced {M}athematics}, 1, 1997.

\bibitem[dJ95]{de_jong:rigid_pi_1}
A.J. de~Jong.
\newblock {\'E}tale fundamental groups of non-archimedean analytic spaces.
\newblock {\em Compositio Mathematica}, 97(1--2):89--118, 1995.

\bibitem[DLLZ19a]{dllz:foundations}
H.~Diao, K.-W. Lan, R.~Liu, and X.~Zhu.
\newblock Logarithmic adic spaces: some foundational results.
\newblock 2019.
\newblock arXiv preprint 1912.09836, version 1 accessed.

\bibitem[DLLZ19b]{dllz:riemann-hilbert}
H.~Diao, K.-W. Lan, R.~Liu, and X.~Zhu.
\newblock Logarithmic {R}iemann--{H}ilbert correspondences for rigid varieties.
\newblock 2019.
\newblock arXiv preprint 1803.05786, version 3 accessed.

\bibitem[DN18]{deglise-niziol}
F.~D{\'e}glise and W.~Nizio{\l}.
\newblock On $p$-adic absolute {H}odge cohomology and syntomic coefficients, i.
\newblock {\em Commentarii Mathematici Helvetici}, 93(1):71--131, 2018.

\bibitem[Fal83]{faltings}
G.~Faltings.
\newblock Endlichkeitss{\"a}tze f{\"u}r abelsche {V}ariet{\"a}ten {\"u}ber
  {Z}ahlk{\"o}rpern.
\newblock {\em Inventiones mathematicae}, 73:349--366, 1983.

\bibitem[Fon94]{fontaine3}
J-M. Fontaine.
\newblock Repr{\'e}sentations $p$-adiques semi-stables.
\newblock {\em Ast{\'e}risque}, 223: P{\'e}riodes $p$-adiques:113--184, 1994.

\bibitem[FvdP04]{fresnel-van_der_put:rigid}
J.\ Fresnel and M.~van~der Put.
\newblock {\em Rigid Analytic Geometry and Its Applications}.
\newblock Number 208 in Progress in Mathematics. Birkh\"auser Basel, 2004.

\bibitem[Had11]{hadian:motivic_pi_1}
M.~Hadian.
\newblock Motivic {F}undamental {G}roups and {I}ntegral {P}oints.
\newblock {\em Duke Mathematical Journal}, 160(3):503--565, 2011.

\bibitem[Hub96]{huber}
R.~Huber.
\newblock {\em {\'E}tale {C}ohomology of {R}igid {A}nalytic {V}arieties and
  {A}dic {S}paces}, volume~30 of {\em Aspects of Mathematics}.
\newblock Springer-Verlag, 1996.

\bibitem[Ill02]{illusie}
L.~Illusie.
\newblock An overview of the work of {K}. {F}ujiwara, {K}. {K}ato and {C}.
  {N}akamura on logarithmic {\'e}tale cohomology.
\newblock {\em Ast{\'e}risque}, 279:271--322, 2002.

\bibitem[Kim05]{minhyong:siegel}
M.~Kim.
\newblock The motivic fundamental group of {$\mathbf
  P^1\setminus\{0,1,\infty\}$} and the theorem of {S}iegel.
\newblock {\em Inventiones Mathematicae}, 161(3):629--656, 2005.

\bibitem[Kim09]{minhyong:selmer}
M.~Kim.
\newblock The unipotent {A}lbanese map and {S}elmer varieties for curves.
\newblock {\em Publications of the Research Institute for Mathematical
  Sciences}, 45:89--133, 2009.

\bibitem[KT08]{minhyong-tamagawa}
M.~Kim and A.~Tamagawa.
\newblock The $l$-component of the unipotent {A}lbanese map.
\newblock {\em Mathematische Annalen}, 340:223--235, 2008.

\bibitem[LZ17]{liu-zhu:riemann-hilbert}
R.~Liu and X.~Zhu.
\newblock Rigidity and a {R}iemann--{H}ilbert correspondence for $p$-adic local
  systems.
\newblock {\em Inventiones Mathematicae}, 207(1):291--343, 2017.

\bibitem[Mil17]{milne}
J.S. Milne.
\newblock {\em Algebraic Groups}.
\newblock Number 170 in Cambridge studies in advanced mathematics. Cambridge
  University Press, 2017.

\bibitem[Ols11]{olsson}
M.C. Olsson.
\newblock Towards non-abelian $p$-adic {Hodge} theory in the good reduction
  case.
\newblock {\em Memoirs of the American Mathematical Society}, 210(990), 2011.

\bibitem[Ols16]{olsson_again}
M.~Olsson.
\newblock The {B}ar {C}onstruction and {A}ffine {S}tacks.
\newblock {\em Communications in Algebra}, 44(7):3088--3121, 2016.

\bibitem[Sch13]{scholze:relative}
P.~Scholze.
\newblock $p$-adic {H}odge theory for rigid-analytic varieties.
\newblock {\em Forum of Mathematics, Pi}, 1:e1, 2013.

\bibitem[Sch16]{scholze:erratum}
P.~Scholze.
\newblock $p$-adic {H}odge theory for rigid-analytic varieties -- corrigendum.
\newblock {\em Forum of Mathematics, Pi}, 4:e6, 2016.

\bibitem[Shi00]{shiho_1}
A.~Shiho.
\newblock Crystalline fundamental groups {I} -- isocrystals on log crystalline
  site and log convergent site.
\newblock {\em Journal of Mathematical Sciences, the University of Tokyo},
  7:509--656, 2000.

\bibitem[Shi20]{shimizu:p-adic_monodromy}
K.~Shimizu.
\newblock A $p$-adic monodromy theorem for de {R}ham local systems, 2020.
\newblock arXiv preprint (arXiv:2003.10951v2).

\bibitem[Sza09]{szamuely:fundamental_groups}
T.\ Szamuely.
\newblock {\em Galois groups and fundamental groups}.
\newblock Number 117 in Cambridge {S}tudies in {A}dvanced {M}athematics.
  Cambridge University Press, 2009.

\bibitem[Vol03]{vologodsky}
V.~Vologodsky.
\newblock Hodge structure on the fundamental group and its application to
  $p$-adic integration.
\newblock {\em Moscow Mathematical Journal}, 3(1):205--247, 2003.

\end{thebibliography}
\bibliographystyle{alpha}

\end{document}